\documentclass[twoside, 12pt]{article}
\usepackage[latin1]{inputenc}
\usepackage{amscd}
\usepackage{amsmath}
\usepackage{amsfonts}
\usepackage{amssymb}
\usepackage{amsthm}
\usepackage{comment}

\renewcommand{\tilde}{\widetilde}
\renewcommand{\hat}{\widehat}

\newcommand{\norm}[1]{\left\| #1 \right\|}
    \DeclareMathOperator{\Ree}{Re}

    \DeclareMathOperator{\Hom}{Hom}

    \DeclareMathOperator{\Cor}{Cor}

\newcommand{\de}{{\rm d}}
\newcommand{\dd}{\, {\rm d}}
\renewcommand{\d}{\displaystyle}
\renewcommand{\geq}{\geqslant}
\renewcommand{\leq}{\leqslant}

\newcommand{\N}{\mathbb{N}}

\newcommand{\R}{\mathbb{R}}

\newcommand{\C}{\mathbb{C}}

\renewcommand{\phi}{\varphi}
\renewcommand{\epsilon}{\varepsilon}
\newcommand{\tq}{\ |\ }
\newcommand{\D}{\mathbb{D}}
\newcommand{\boL}{\mathcal{L}}

\newcommand{\boB}{\mathcal{B}}

\newtheorem{thm}{Theorem}[section]
\newtheorem*{thm*}{Theorem}
\newtheorem{prop}[thm]{Proposition}
\newtheorem{lem}[thm]{Lemma}
\newtheorem{cor}[thm]{Corollary}

\theoremstyle{definition}

\title{Sharp polynomial estimates for the decay of correlations}
\author{Sébastien Gouëzel\footnote{Université d'Orsay, e-mail
    \texttt{Sebastien.Gouezel@math.u-psud.fr}}}
\date{February 13, 2002}
\begin{document}

\maketitle

\begin{abstract}
We generalize a method developed by Sarig to obtain polynomial
lower bounds for correlation functions for maps with a countable
Markov partition. A consequence is that LS Young's estimates on
towers are always optimal. Moreover, we show that, for functions
with zero average, the decay rate is better, gaining a factor
$1/n$. This implies a Central Limit Theorem in contexts where it
was not expected, e.g. $x+Cx^{1+\alpha}$ with $1/2\leq \alpha<1$.
The method is based on a general result on renewal sequences of
operator, and gives an asymptotic estimate up to any precision of
such operators.
\end{abstract}

\section{Statement of results}

In recent years, several methods have been developed to obtain
polynomial upper bounds for the correlations of some dynamical
systems. However, there was no general method to get polynomial
lower bounds for the decay of correlations, until
Omri Sarig's recent article \cite{sarig:decay}. He used an abstract
result on renewal sequences of operators to obtain lower bounds on the
decay of correlations for Markov maps. As an application, he
proved that the upper bounds obtained by Young on tower maps in
\cite{lsyoung:recurrence} are in many cases optimal. The
goal of this article is to remove some unnecessary assumptions in
\cite{sarig:decay}, and as a
consequence to prove that Young's estimates are optimal in full
generality.

In this article, $\D$ will always denote $\{z\in \C \tq |z|<1\}$.
The analogue of Sarig's theorem on renewal sequences that we
obtain is the following:
\begin{thm}
\label{theoreme_renouvellement}
Let $T_n$ be bounded operators on a Banach space $\boL$ such that
$T(z)=I+\sum_{n\geq 1}z^nT_n$ converges in $\Hom(\boL,\boL)$ for every
$z\in \D$. Assume that:
\begin{enumerate}
\item \textbf{Renewal equation:}
for every $z\in \D$, $T(z)=(I-R(z))^{-1}$ where $R(z)=\sum_{n\geq
  1}z^nR_n$, $R_n\in \Hom(\boL,\boL)$ and $\sum \norm{R_n}<+\infty$.
\item \textbf{Spectral gap:}
$1$ is a simple isolated eigenvalue of $R(1)$.
\item \textbf{Aperiodicity:}
for every $z\in \overline{\D}-\{1\}$, $I-R(z)$ is invertible.
\end{enumerate}
Let $P$ be the eigenprojection of $R(1)$ at $1$. If
$\sum_{k>n}\norm{R_k}=O(1/n^\beta)$ for some $\beta>1$ and $PR'(1)P\not=0$,
then for all n
    \begin{equation*}
    T_n=\frac{1}{\mu}P+\frac{1}{\mu²}\sum_{k=n+1}^{+\infty}P_k+E_n
    \end{equation*}
where $\mu$ is given by $PR'(1)P=\mu P$, $P_n=\sum_{l>n}PR_lP$ and $E_n\in
\Hom(\boL,\boL)$ satisfy
    \begin{equation*}
    \norm{E_n}=\left\{\begin{array}{ll}
    O(1/n^\beta)&\text{if }\beta>2\\
    O(\log n/n^2)&\text{if }\beta=2\\
    O(1/n^{2\beta-2})&\text{if }2>\beta>1
    \end{array}\right.
    \end{equation*}
\end{thm}
Note that, in all cases, $\norm{E_n}=o(1/n^{\beta-1})$, which is
what is needed to obtain sharp asymptotics for the decay of
correlations.
This theorem extends Sarig's: he assumed $\beta>2$
while we only need $\beta>1$. Moreover, the result we obtain is
slightly stronger than Sarig's even in the case $\beta>2$ because
the error term is a $O(1/n^\beta)$ instead of a $O(1/n^{\lfloor
\beta \rfloor})$.

Finally, our aperiodicity assumption
is weaker than Sarig's who needed to suppose that the spectral
radius of $R(z)$ were $<1$ for every $z\not=1$. Our assumption is
necessary because other eigenvalues equal to $1$ would generate
other terms in the asymptotic expression of $T_n$ (which could be
calculated using the same methods as in the following proof, and
would involve the spectral projection at these points). For example,
if $R(z)=z^2$, then $T_{2n}=1$ while $T_{2n+1}=0$, which shows that
the conclusions of the theorem are not valid any more (there is a
periodicity problem).
This less
restrictive aperiodicity hypothesis will be useful for example
when applied to tower maps (see Corollary
\ref{cor_LSYoung_tower}).

It is in fact possible to give an asymptotic estimate of $T_n$ up
to an error term $O(1/n^\beta)$ even when $\beta \leq 2$.
However, the result is quite technical to state, and will be
deferred to Section \ref{section_termes_suivants}. The following
consequence of Theorem \ref{Tn_ordre_plus_grand} will be
sufficient for most dynamical applications.

\begin{thm}
\label{thm_quand_Pf_nul}
Under the hypotheses of Theorem \ref{theoreme_renouvellement}, if
$f\in \boL$ is such that $Pf=0$, then $\norm{T_n f}=O(1/n^\beta)$.
\end{thm}

These abstract results enable us to enhance the applications in
\cite{sarig:decay}. We state briefly the results we obtain,
without recalling all the notations. In Section \ref{markov_maps},
a precise meaning will be given to all the notions involved. The
following theorem is stated more precisely as Theorem
\ref{thm_applications_markov_precis}.

\begin{thm}
\label{thm_markov}
Let $(X,\boB,m,T,\alpha)$ be a topologically mixing
 probability preserving Markov
map, and $\gamma\subset\alpha$. Denote by $T_\gamma$ the map
induced by $T$ on $Y=\bigcup \gamma$ -- it is a Markov map for a
subpartition $\delta$ of $\gamma$. Assume that the distortion of
$T_\gamma$ decreases exponentially fast and that $T_\gamma$ has
the ``big image'' property (which is always true when $\gamma$ is
finite). Assume moreover that $m[\phi_\gamma>n]=O(1/n^\beta)$ for
some $\beta>1$, where $\phi_\gamma$ is the first return time from
$Y$ to $Y$.

Then $\exists \theta\in
(0,1)$, $C>0$ such that $\forall f,g$ integrable and supported inside $Y$,
    \begin{equation*}
    \left|\Cor(f,g\circ T^n)-\left(
    \sum_{k=n+1}^{\infty}m[\phi_\gamma>k]\right)\int f \int
    g \right| \leq C F_\beta(n) \norm{g}_\infty \norm{f}_{\boL}
    \end{equation*}
where $F_\beta(n)=1/n^\beta$ if $\beta>2$, $\log n/n^2$ if $\beta=2$ and
$1/n^{2\beta-2}$ if $2>\beta>1$ (and $\boL$ denotes the space of
$\theta$-Hölderian functions on $Y$).

Moreover, if $\int f=0$, then $\Cor(f,g\circ T^n)=O(1/n^\beta)$.
\end{thm}

When $m[\phi_\gamma>k]\asymp 1/n^\beta$ and $\int f, \int
g\not=0$, Theorem \ref{thm_markov} implies that $\Cor(f,g\circ
T^n)\asymp 1/n^{\beta-1}$. Thus, the exact speed of decay of
correlations is polynomial, with exponent $\beta-1$.
Surprisingly, the decay rate is better for functions with zero
integral, with a gain of $1$ in the exponent. This kind of result
is to the knowledge of the author new, and does not seem to be
obtainable by more crude estimates: the methods giving only upper
bounds on the speed of decay of correlations do not distinguish
between functions with zero or non-zero integral, since they do
not ``see'' the higher order terms in the expansion of $T_n$.

As an application, we obtain the summability of the correlations
for functions with zero integral (and supported in $Y$) even when
$\beta \leq 2$, which gives a Central Limit Theorem in cases where it was not
expected. Note that the condition of zero integral is important
and can not be eliminated by subtracting a constant, since the
functions would not remain supported in $Y$. In fact, the
estimate in the previous theorem shows that, when $\beta\leq 2$,
the correlations are \emph{not} summable for a function with
nonzero integral supported in $Y$, which makes it very unlikely
for a CLT to hold (and replacing $f$ by $f-\frac{\int
f}{m(Y)}1_Y$ does not give any result since there is no good
control on the behavior of $1_Y$). In the same way, this speed of
decay of correlations does not hold for general functions with
zero integral but not supported in $Y$: take a function $f$ of
nonzero integral supported in $Y$, the function $g:=f-\int f$ has
zero integral but its correlations are the same as those of $f$,
whence they decay at a rate $\asymp 1/n^{\beta-1}$.

The following Central Limit Theorem is stated more accurately as
Theorem \ref{markov_CLT}.
\begin{thm}
Under the same hypotheses as in Theorem \ref{thm_markov}, if
$f\in \boL$ is supported in $Y$ and $\int f=0$, then the sequence
$\frac{1}{\sqrt{n}} \sum_{k=0}^{n-1} f\circ T^k$ converges in law
to a Gaussian random variable of zero mean and finite variance
$\sigma$, with
  \begin{equation*}
  \sigma^2=-\int f^2\dd m+2\sum_{n=0}^\infty \int f\cdot f\circ
  T^n \dd m
  \end{equation*}
\end{thm}

Finally, even though Theorem \ref{thm_markov} describes the speed
of decay of correlations only for functions $f$ and $g$ supported
in $Y$, it is possible to drop this hypothesis on $g$. However,
the results obtained are less precise and give only an upper
bound on the decay of correlations, in $O(1/n^{\beta-1})$ if
$\int f\not=0$ and in $O(1/n^\beta)$ if $\int f=0$ (see Theorem
\ref{decay_sur_tout_X} and Proposition
\ref{decay_sur_tout_X_integrale_nulle}). This kind of results is
useful in the proof of the Central Limit Theorem.

The following corollaries are already present in weaker form in
\cite{sarig:decay}, where the notations are explained. Some
details on their proofs will be given in the last section of this
article. The first corollary (stated more precisely as Corollary
\ref{cor_LSV}) deals with an explicit one-dimensional Markov map
with a neutral fixed point, while the second corollary (see
section \ref{section_LSY} and Corollary \ref{cor_LSY}) is
essentially Theorem \ref{thm_markov} expressed in the framework
of LS Young towers, which are devices built up from non-Markov
maps which have proved very useful in studying their statistical
properties (see \cite{lsyoung:recurrence}).

\begin{cor}
\label{cor_liv}
In the case of the Liverani-Saussol-Vaienti map $T:[0,1]\to
[0,1]$ defined by
   \begin{equation*}
   T(x)=\left\{ \begin{array}{cl}
   x(1+2^\alpha x^\alpha) &\text{if }0\leq x\leq 1/2
   \\
   2x-1 &\text{if }1/2<x\leq 1
   \end{array}\right.
   \end{equation*}
(see \cite{liverani_saussol_vaienti}), if $\alpha\in (0,1)$, $f$
is Lipschitz, $g$ is bounded measurable, $\int f,\int g\not=0$ and
$f,g=0$ in a neighborhood of $0$, then $\Cor(f,g\circ T^n)\sim
\frac{1}{4}h(\frac{1}{2}) \alpha^{-1/\alpha}
(\frac{1}{\alpha}-1)^{-1}n^{1-1/\alpha}\int f\int g$ with respect
to the invariant probability measure.

Moreover, if $\int f=0$ (and $f,g$ are still Lipschitzian and
zero in a neighborhood of $0$), then $\Cor(f,g\circ
T^n)=O(1/n^{1/\alpha})$. Consequently, $f$ satisfies a Central
Limit Theorem.
\end{cor}
This result is in fact not specific to this particular map and can
easily be extended to a class of maps admitting a neutral fixed
point in $0$ with a prescribed behavior, and expanding outside of
any neighborhood of $0$, making use of the following corollary
and the techniques of \cite{lsyoung:recurrence}.

\begin{cor}
\label{cor_LSYoung_tower}
Let $(\Delta,\boB,m,F)$ be a probability preserving LS Young
tower with $\gcd\{R_i\}=1$ and $m[R>n]=O(1/n^\beta)$ where
$\beta>1$. If $f\in C_{\theta}(\Delta)$, $g\in L^{\infty}$ are
supported inside $\bigcup_0^{N-1}\Delta_i$ for some $N$, then
$\Cor(f,g\circ F^n)=\sum_{k>n} m[R>k]\int f\int
g+O(F_\beta(n))$.

Moreover, if $\int f=0$, then $\Cor(f,g\circ F^n)=O(1/n^\beta)$.
Thus, $f$ satisfies a Central Limit Theorem.
\end{cor}
The aperiodicity hypothesis on $\gcd\{R_i\}$ is the same as
Young's, and cannot be omitted. In her paper
\cite{lsyoung:recurrence}, Young proved that, if
$m[R>n]=O(1/n^{\beta-1})$, then $\Cor(f,g\circ
F^n)=O(1/n^{\beta-1})$ for any $f\in C_{\theta}(\Delta)$, $g\in
L^{\infty}$ (not necessarily supported in
$\bigcup_0^{N-1}\Delta_i$). Corollary \ref{cor_LSYoung_tower}
proves that this upper estimate is in fact optimal, and gives
additionally a Central Limit Theorem even if $1<\beta\leq 2$.

From this point on, the paper is divided into two parts: the
first one (sections \ref{section_prelim},
\ref{preuve_main_lemma}, \ref{section_preuve_main} and
\ref{section_termes_suivants}) is devoted to the proof of the
abstract results on renewal sequences of operators, and the
second one (sections \ref{markov_maps} and
\ref{section_applications}) deals with the application to Markov
maps.

\section{Preliminary results}
\label{section_prelim}

\subsection{$C^{1+\alpha}$ functions in Banach algebras}
The results in this section are mainly straightforward
computations, and most of them can be found in \cite{sarig:decay}.

Let $\boB$ be a Banach algebra (in the applications of interest to
us, $\boB=\Hom(\boL,\boL)$). Fix $K$ a compact subset of $\C$.
The distance on $K$ will not be the usual one, but the geodesic
distance, i.e. $d(x,y)$ is the infimum of the lengths of
$C^1$-paths in $K$ joining $x$ to $y$. We assume that this
distance is equivalent to the usual one, which will be true for
$K=\overline{\D}$ or $K=S^1$.

Fix some $0<\alpha<1$.
For any $f:K\to\boB$, we will say that $f$ is $C^\alpha$ if there exists a
constant $C$ such that, for any $x,y\in K$, $\norm{f(x)-f(y)}\leq C
d(x,y)^\alpha$. Let $D_\alpha(f)$ denote the least such constant. We write
$\norm{f}_\alpha=\norm{f}_\infty+D_\alpha(f)$, and denote by
$C^\alpha(K)$ the space of all functions such that $\norm{f}_\alpha<+\infty$.

\begin{prop}
The space $(C^\alpha(K),\norm{\ }_\alpha)$ is a Banach algebra. In fact, we
even have, for $f,g\in C^\alpha(K)$, $D_\alpha(fg)\leq \norm{f}_\infty
D_\alpha(g)+\norm{g}_\infty D_\alpha(f)$.
\end{prop}

We say that $f:K\to \boB$ is $C^1$ if there exists a continuous
function $g:K\to\boB$ such that $f(x+h)-f(x)-h g(x)=o(h)$ for any
$x \in K$. The function $g$ is unique if it exists, and we write
$g=f'$.

\begin{prop}
If $f$ is $C^1$ on $K$, then $D_\alpha(f)\leq
\norm{f'}_\infty$.
\end{prop}
\begin{proof}
Let $x,y\in K$ with $d(x,y)<1$. Let $\gamma$ be a $C^1$ path in $K$ from
$x$ to $y$. The Taylor-Lagrange inequality along this path gives
$\norm{f(x)-f(y)}\leq \norm{f'}_\infty l(\gamma)$.
\end{proof}
We consider the geodesic distance on $K$ instead of the usual one
precisely to get the above proposition.

Let $C^{1+\alpha}(K)$ denote the space of all $C^1$ functions from $K$ to
$\boB$ whose derivative is $C^\alpha$, endowed with the norm
$\norm{f}_{1+\alpha}=\norm{f}_\infty+\norm{f'}_\infty+\frac{1}{2}
D_\alpha(f')$.

\begin{prop}
The space $(C^{1+\alpha}(K),\norm{\ }_{1+\alpha})$ is a Banach algebra.
\end{prop}

The following proposition will be used systematically in Section
\ref{preuve_main_lemma}, often without explicit reference.
\begin{prop}
\label{borne_norme_1alpha_inverse}
Let $f:K\to\boB$ be a $C^{1+\alpha}$ function such that, for every $z\in K$,
$f(z)$ is invertible (as an element of $\boB$). If $g(z)=f(z)^{-1}$,
then $g$ is $C^{1+\alpha}$ and there is an inequality $\norm{g}_{1+\alpha}\leq
F(\norm{g}_\infty,\norm{f}_{1+\alpha})$ for some universal polynomial function $F$.
\end{prop}
\begin{proof}
Differentiating $g(z)=f(z)^{-1}$, we get $g'(z)=-g(z) f'(z)g(z)$,
hence $\norm{g'}_\infty \leq \norm{g}_\infty^2 \norm{f'}_\infty$.

Then we note that $\norm{g'}_\alpha \leq \norm{g}_\alpha \norm{f'}_\alpha
\norm{g}_\alpha \leq \norm{g'}_\infty \norm{f'}_\alpha \norm{g'}_\infty$.
The control on $\norm{g'}_\infty$ enables us to conclude.
\end{proof}

\subsection{Fourier series in Banach algebras}

Let $\boB$ be a Banach algebra. For $f:S^1\to \boB$ a continuous
function, we define the $n^{\text{th}}$ Fourier coefficient of $f$ to
be the
element of $\boB$ defined by
    \begin{equation*}
    c_n(f)=\frac{1}{2\pi}\int_0^{2\pi} f(e^{i\theta})e^{-in\theta}d\theta
    \end{equation*}

Let us first recall a very useful result concerning functions from
$S^1$ to $\C$.
\begin{thm}[Wiener Lemma]
Let $f:S^1\to\C$ be a continuous function, everywhere nonzero,
whose Fourier coefficients are summable. Then the Fourier coefficients
of $1/f$ are also summable.
\end{thm}
The proof of this result, which uses commutative Banach algebra
techniques, can be found for example in \cite{katznelson}.

\begin{prop}
If $f:S^1\to\boB$ is continuous and satisfies $\sum
\norm{c_n(f)}<+\infty$, then $f(e^{i\theta})=\sum
c_n(f)e^{in\theta}$, the series converging in norm.
\end{prop}
\begin{proof}
Replacing $f$ by $f-\sum c_n(f) e^{in\theta}$, we can assume that $c_n(f)=0$
for every $n$, and we want to prove that $f=0$.

Suppose on the contrary the existence of $z$ such that
$f(z)\not=0$. There exists a linear form $\phi$ on $\boB$ with
$\phi(f(z))\not=0$. The linearity of $\phi$ gives $c_n(\phi \circ
f)=\phi(c_n(f))=0$ for every $n$. As $\phi\circ f$ is
complex-valued, a classical result (proved for example using
Parseval's equality) gives $\phi\circ f=0$, which is a
contradiction.
\end{proof}

\begin{prop}
\label{stable_produit}
If the Fourier coefficients of $f$ and $g$
are summable, then it is also the case of $fg$.
\end{prop}
\begin{proof}
Write $f=\sum c_n(f)e^{in\theta}$ and $g=\sum c_n(g)e^{in\theta}$.
Multiply, integrate (and permute: everything converges in norm) to
get $c_n(fg)=\sum_{k+l=n}c_k(f)c_l(g)$. As a consequence,
    \begin{equation*}
    \sum \norm{c_n(fg)}\leq \sum_{k,l}\norm{c_k(f)c_l(g)}\leq
    \sum_k \norm{c_k(f)}\sum_l\norm{c_l(g)}
    \end{equation*}
\end{proof}

\begin{prop}
\label{C_1alpha_donc_coeffs_sommables}
If $f:S^1\to \boB$ is
$C^{1+\alpha}$ then $\norm{c_n(f)}\leq C
\frac{\norm{f}_{1+\alpha}}{n^{1+\alpha}}$ for some universal
constant $C$.
\end{prop}
The classical proof for complex valued functions can be found in
\cite{katznelson} and is easily adapted to this context (see also
\cite[Lemma 3]{sarig:decay}).

\section{Proof of Sarig's first main lemma under our weaker
  assumptions
\label{preuve_main_lemma}}
The following lemma, which is the analogue of Sarig's first main
lemma,
is crucial to the proof of Theorem
\ref{theoreme_renouvellement}.

\begin{lem}[First Main Lemma]
\label{first_main_lemma}
Under the assumptions of Theorem
\ref{theoreme_renouvellement},
    \begin{equation*}
    \sum_{n=1}^{\infty}\norm{T_n-T_{n+1}}<\infty
    \end{equation*}
\end{lem}
As $(I-R(z))^{-1}=\sum T_n z^n$ (where we write $T_0=I$), we have
$A(z):=(1-z)(I-R(z))^{-1}=\sum (T_n-T_{n-1})z^n$. Our strategy is
to study $A$ on $S^1$, and to see that its Fourier coefficients
are summable. As $A(z)=\left(\frac{I-R(z)}{1-z}\right)^{-1}$, $A$
is well controlled on $S^1$ outside of any neighborhood of $1$.
Near $1$, the problem comes from the eigenvalue $\lambda(z)$ of
$R(z)$ closest to $1$. To use Fourier series methods to control
this eigenvalue, we must be able to extend $\lambda(z)$ to the
whole circle $S^1$; that is why we will have to modify $R(z)$ and
to construct a function $\tilde{R}(z)$ on $S^1$, whose spectrum
will be ``nice.''

\emph{Proof of Lemma \ref{first_main_lemma}.}

We will write $\beta=1+\alpha$. We can assume $0<\alpha<1$, which
amounts only to weakening the hypotheses.

\emph{Step 1: $R(z)$ is $C^{1+\alpha}$ on $\overline{\D}$.}
\begin{proof}
As $R(z)=\sum R_n z^n$ with $\sum_{k>n}
\norm{R_k}=O(1/n^{1+\alpha})$, we have $R_n=O(1/n^{1+\alpha})$,
and the series defining $R$ converges in norm on all
$\overline{\D}$. Thus, $R$ is continuous on $\overline{\D}$.

The sum $F(z)=\sum  n R_n z^{n-1}$ converges also in norm on
$\overline{\D}$, as $\sum_n n\norm{R_n}=\sum_{n\geq 1}
S_n<+\infty$ (where we write $S_n=\sum_{k\geq n}\norm{R_k}$).
Hence, this is the derivative of $R$ and $R$ is $C^1$ on
$\overline{\D}$.

What remains to be checked is that $F$ is $C^\alpha$. Let $z$ and
$z+h$ be two points in $\overline{\D}$, we estimate from above
$\norm{F(z+h)-F(z)}$. The Taylor-Lagrange inequality gives, for
every $n\in \N$, $|(z+h)^n-z^n|\leq n|h|$. Let $N\in \N$. We have
    \begin{align*}
    |F(z+h)-F(z)|&
    \leq \sum_{n=0}^N n \norm{R_n} |(z+h)^{n-1} -z^{n-1}|
    +\sum_{n=N+1}^{+\infty} 2 n \norm{R_n}
    \\&
    \leq \sum_{n=0}^N n^2(S_n-S_{n+1})|h|
    +\sum_{n=N+1}^{+\infty} 2n (S_n-S_{n+1})
    \\&
    \leq \sum_{n=0}^{N-1}2n S_n |h|
    +\sum_{n=N+2}^{+\infty}2S_n +2(N+1) S_{N+1}
    \end{align*}
As $nS_n=O(1/n^\alpha)$,
$\sum_{n=0}^{N-1}nS_n=O(1/N^{\alpha-1})$, while
$\sum_{n=N+1}^{+\infty} S_n=O(1/N^\alpha)$ and
$NS_N=O(1/N^\alpha)$. Hence, for some constants $C$ and $D$
(independent of $N$ or $h$),
    \begin{equation*}
    |F(z+h)-F(z)|
    \leq \frac{C}{N^{\alpha-1}}|h|+\frac{D}{N^\alpha}
    \end{equation*}
If we choose $N$ close to $1/|h|$, we get a bound of the order of
$|h|^\alpha$.
\end{proof}

\emph{Step 2: $\d \frac{R(z)-R(1)}{z-1}$ can be continuously extended
  to $S^1$, and its Fourier coefficients are summable.}
\begin{proof}
For $z\not=1$,
    \begin{equation*}
    \frac{R(z)}{1-z}=\frac{\sum R_n
    z^n}{1-z}=\sum_{n=0}^{+\infty}\left(\sum_{k=0}^n R_k\right) z^n
    \end{equation*}
Moreover, $R(1)/(1-z)=\sum_{n=0}^{+\infty}(\sum_{k=0}^{+\infty}
R_k) z^n$, hence
    \begin{equation*}
    \frac{R(z)-R(1)}{z-1}=\sum_{n=0}^{+\infty}\left(\sum_{k=n+1}^{+\infty}
    R_k\right) z^n
    \end{equation*}
The last sum converges in norm, because
$\sum_{k=n+1}^{+\infty}\norm{R_k}=O(1/n^{1+\alpha})$ is summable.
This guarantees a continuous extension to $1$. Moreover, the
$n^{\text{th}}$ Fourier coefficient is
$\sum_{k=n+1}^{+\infty}R_k$, which is summable.
\end{proof}

\emph{Step 3: Construction of a function $\tilde{R}$ on $S^1$,
  equal to $R$ in a neighborhood of $1$, $C^{1+\alpha}$ outside of any
  neighborhood of $1$, whose spectrum consists in an isolated eigenvalue
  $\tilde{\lambda}(z)$ close to $1$ together with a compact subset of
  $\C-\{1\}$, with $\tilde{\lambda}(z)\not=1$ for $z\not=1$.
  Furthermore, for any $\epsilon>0$,
  $\tilde{R}$ can be chosen such that $\forall z\in S^1$,
  $\|\tilde{R}(z)-R(1)\|<\epsilon$.}

\begin{proof}
We construct two candidates for $\tilde{R}$, $U$ and $V$. The
second one, i.e. $V$, will be the good one.

Fix some $\gamma>0$, very small. Let $\phi+\psi$ be a $C^\infty$
partition of unity associated to the sets $\{\theta \in
[0,\gamma)\}$ and $\{\theta \in (\gamma-\eta, \pi/2]\}$ where
$\theta$ is the angle on the circle (for some very small
$0<\eta<\gamma$). We define $U(z) =\phi(z)R(z)
+\psi(z)R(e^{i\gamma})$ on $\{\theta\in [0,\pi/2]\}$ : $U$ is
equal to $R$ on $\{\theta\in [0,\gamma-\eta]\}$ and to
$R(e^{i\gamma})$ on $\{\theta\in [\gamma,\pi/2]\}$. In particular,
the spectrum of $U(z)$ will be "almost the same" as the spectrum
of $R(1)$, if $\gamma$ is small enough.

We define in the same way $U$ on $\{\theta\in [-\pi/2,0]\}$, equal
to $R(e^{-i \gamma})$ on $\{\theta\in [-\pi/2,-\gamma]\}$ and to
$R$ on $\{\theta\in [-\gamma+\eta,0]\}$.

Finally, we construct $U$ on the remaining half-circle by
symmetrizing, i.e. $U(e^{i(\pi/2+a)})=U(e^{i(\pi/2-a)})$, to
ensure that everything fits well.

There is a well defined eigenvalue close to $1$ for every $U(z)$,
depending continuously on $z$, which we denote by $\rho(z)$. The
problem would be solved if $\rho(z)\not=1$ for $z\not=1$, which
is not the case since $\rho(-1)=\rho(1)=1$. Consequently, we have
to perturb $\rho$ a little. There exists a $C^\infty$ function
$\nu$ on $\{\theta\in [\pi/2,3\pi/2]\}$ arbitrarily close to
$\rho$. We can assume that $\nu$ is transversal to $\{1\}$, i.e.
that it does not take the value $1$. On $\{\theta\in
[\pi/2+\eta,3\pi/2-\eta] \}$, we define
$V(z)=\frac{\nu(z)}{\rho(z)}U(z)$ : its
  eigenvalue close to $1$ is $\nu(z)\not=1$. Finally, we glue
  $U$ and $V$ together on $\{\theta\in [\pi/2,\pi/2+\eta]\}$ and
$\{\theta\in [3\pi/2-\eta,3\pi/2]\}$ with a partition of unity, as
above. As the spectrum of $U(e^{i\pi/2})=R(e^{i\gamma})$ does not
contain $1$, the gluing will not give an eigenvalue equal to $1$
if we choose $\eta$ small enough and $\nu$ close enough to $\rho$.
\end{proof}

\emph{Step 4: $\d \frac{\tilde{R}(z)-\tilde{R}(1)}{z-1}$
can be continuously extended to $S^1$ and its Fourier coefficients are
summable.}
\begin{proof}
As $\tilde{R}(1)=R(1)$,
    \begin{equation*}
    \frac{\tilde{R}(z)-\tilde{R}(1)}{z-1}
    =\frac{\tilde{R}(z)-R(z)}{z-1} + \frac{R(z)-R(1)}{z-1}
    \end{equation*}
The first term is $C^{1+\alpha}$ outside of any neighborhood of
$1$, and zero on a neighborhood of $1$. Thus, it is
$C^{1+\alpha}$, which shows that its Fourier coefficients are
summable by Proposition
 \ref{C_1alpha_donc_coeffs_sommables}.

The coefficients of the second term $\frac{R(z)-R(1)}{z-1}$ are
summable by Step 2, which gives the conclusion.
\end{proof}

\emph{Step 5: Let $\tilde{P}(z)$ denote the spectral projection of
  $\tilde{R}(z)$ corresponding to its eigenvalue $\tilde{\lambda}(z)$
  close to $1$.Then $\tilde{P}(z)$ is $C^{1+\alpha}$, and its Fourier
  coefficients are summable.}
\begin{proof}
The projection $\tilde{P}(z)$ can be written, for $\delta$ small
enough (and independent of $z$ if, in Step 3, $\epsilon$ was taken
small enough),
    \begin{equation*}
    \tilde{P}(z)=\frac{1}{2i\pi}\int_{|u-1|=\delta}
    \frac{1}{uI-\tilde{R}(z)}\dd u
    \end{equation*}
We already know that $\tilde{R}$ is $C^{1+\alpha}$, which is also
true of $uI-\tilde{R}$ for every $u$, and of
$(uI-\tilde{R})^{-1}$ (with a uniform bound on its $C^{1+\alpha}$
norm) by Proposition \ref{borne_norme_1alpha_inverse}. So, we can
integrate to get a $C^{1+\alpha}$ function.

The summability of the coefficients is then a corollary of
Proposition \ref{C_1alpha_donc_coeffs_sommables}.
\end{proof}

\emph{Step 6: The function $\d
\frac{\tilde{P}(z)-\tilde{P}(1)}{z-1}$ can be continuously
extended to $S^1$ and its Fourier coefficients are summable.}

\begin{proof}
The expression of the spectral projection used in Step 5 gives,
after integration,
    \begin{align*}
    \frac{\tilde{P}(z)-\tilde{P}(1)}{z-1}&
    =\frac{1}{2i\pi}\int_{|u-1|=\delta}
    \frac{1}{uI-\tilde{R}(z)}\frac{\tilde{R}(z)
    -\tilde{R}(1)}{z-1}\frac{1}{uI-\tilde{R}(1)}\dd u
    \end{align*}

Let us fix $u$ such that $|u-1|=\delta$. We have seen in Step 5
that the coefficients of $\frac{1}{uI-\tilde{R}(z)}$ were
summable. Moreover, Step 4 gives the summability of the
coefficients of $\d \frac{\tilde{R}(z)-\tilde{R}(1)}{z-1}$. As a
consequence, the coefficients of the product
$\frac{1}{uI-\tilde{R}(z)}\frac{\tilde{R}(z)
    -\tilde{R}(1)}{z-1}$ are also summable.

To obtain the summability of the coefficients of  $\d
\frac{\tilde{P}(z)-\tilde{P}(1)}{z-1}$, we just have to integrate
with respect to $u$, since
    \begin{equation*}
    c_n\left(\frac{\tilde{P}(z)-\tilde{P}(1)}{z-1}\right)
    =\frac{1}{2i\pi}\int_{|u-1|=\delta}
    c_n \left(\frac{1}{uI-\tilde{R}(z)}\frac{\tilde{R}(z)
    -\tilde{R}(1)}{z-1}\frac{1}{uI-\tilde{R}(1)}\right)\dd u
    \end{equation*}
To conclude, we must get a uniform summable bound on the Fourier
coefficients in the integral, i.e. we have to check that all
previous estimates are uniform in $u$, which does not present any
difficulty: the norms of $(uI-\tilde{R}(z))^{-1}$, for
$|u-1|=\delta$ and $z\in S^1$, are bounded by compactness, and so are the
$1+\alpha$ norms of $uI-\tilde{R}(z)$. Proposition
\ref{borne_norme_1alpha_inverse} guarantees that the $1+\alpha$
norms of $(uI-\tilde{R})^{-1}$ are bounded by a constant
independent of $u$. Proposition
\ref{C_1alpha_donc_coeffs_sommables} gives that
$c_n((uI-\tilde{R})^{-1})=O(1/n^{1+\alpha})$ uniformly in $u$,
which enables us to conclude.
\end{proof}

\emph{Step 7: $\d \frac{\tilde{\lambda}(z)-1}{z-1}\to \mu$ as
$z\to 1$ on $S^1$, where $\mu\not=0$ is such that
$P(1)R'(1)P(1)=\mu P(1)$. Hence, the function $\d
\frac{z-1}{\tilde{\lambda}(z)-1}$ is well defined. Moreover, its
Fourier coefficients are summable.}

\begin{proof}
For every $z\in S^1-\{1\}$, we have
    \begin{equation}
    \label{exprime_valeur_propre}
    \frac{1-\tilde{\lambda}(z)}{1-z}\tilde{P}(z)
    =\frac{I-\tilde{R}(z)}{1-z}\tilde{P}(z)
    =\frac{\tilde{R}(1)-\tilde{R}(z)}{1-z}\tilde{P}(z)+(I-R(1))
    \frac{\tilde{P}(z)-\tilde{P}(1)}{1-z}
    \end{equation}
If we multiply on the left by $\tilde{P}(z)$ and let $z$ go to
$1$, the righthand term tends to $P(1)R'(1)P(1)$ (because the
other term tends to $P(1)(I-R(1))P'(1)=0$, and we can drop the
tildes because $R=\tilde{R}$ in a neighborhood of $1$). But
$P(1)R'(1)P(1)$ can be written $\mu P(1)$, with $\mu\not=0$
according to the hypotheses. We get
    \begin{equation*}
    \frac{1-\tilde{\lambda}(z)}{1-z}\tilde{P}(z)\xrightarrow[z\to 1]{}
    \mu P(1)
    \end{equation*}
Apply a linear form $\phi$ such that $\phi
(\tilde{P}(z))\not=0$ for every $z$ (which is possible: take
$\phi(P(1))\not=0$, and then $\epsilon$ small enough in the construction of $\tilde{R}$).
We obtain the convergence of $\frac{1-\tilde{\lambda}(z)}{1-z}$ to $\mu$.

Then, we show that the Fourier coefficients of the continuous
function $\frac{1-\tilde{\lambda}(z)}{1-z}$ are summable. In
Equation \eqref{exprime_valeur_propre}, all terms on the righthand
side have their coefficients summable, according to the preceding
steps. This remains true when we apply $\phi$, i.e.
$\frac{1-\tilde{\lambda}(z)}{1-z}\phi(\tilde{P}(z))$ has its
coefficients summable. In the same way, $\tilde{P}(z)$ has its
coefficients summable, and $\phi(\tilde{P}(z))$ too. But this is
a complex function, everywhere nonzero, so the Wiener lemma gives
that its inverse $1/\phi(\tilde{P}(z))$ has also summable
coefficients. Multiplying, we obtain the summability of the
coefficients of $\frac{1-\tilde{\lambda}(z)}{1-z}$.

Using once more the Wiener lemma (since
$\frac{1-\tilde{\lambda}(z)}{1-z}$ is everywhere nonzero by
construction of $\tilde{R}$), we get the conclusion.
\end{proof}

\emph{Step 8: $(z-1)(\tilde{R}(z)-I)^{-1}$ can be continuously
  extended to $1$, and its Fourier coefficients are summable.}
\begin{proof}
Let $\tilde{Q}(z)$ denote the spectral projection
$I-\tilde{P}(z)$. Then, for every $z\not=1$,
    \begin{equation}
    \label{prolongement_en_1}
    \begin{split}
    (1-z)(I-\tilde{R}(z))^{-1}
    &=\frac{1-z}{1-\tilde{\lambda}(z)}\tilde{P}(z)
    +(1-z)(I-\tilde{R}(z))^{-1}\tilde{Q}(z)
    \\&=\frac{1-z}{1-\tilde{\lambda}(z)}\tilde{P}(z) +
    (1-z)(I-\tilde{R}(z)\tilde{Q}(z))^{-1} \tilde{Q}(z)
    \end{split}
    \end{equation}
$I-\tilde{R}(z)\tilde{Q}(z)$ is everywhere invertible on $S^1$
and is $C^{1+\alpha}$ (this is true for $\tilde{Q}$ because
$\tilde{P}$ is $C^{1+\alpha}$ by Step 5 and
$\tilde{Q}=I-\tilde{P}$). Proposition
\ref{borne_norme_1alpha_inverse} gives that its inverse is
$C^{1+\alpha}$, hence its coefficients are summable, which
remains true when it is multiplied by $\tilde{Q}(z)$ which is
$C^{1+\alpha}$.

To conclude, we have to show that
$\frac{1-z}{1-\tilde{\lambda}(z)}\tilde{P}(z)$ has its Fourier
coefficients summable. We already know this for $\tilde{P}(z)$
(Step 5) and $\frac{1-z}{1-\tilde{\lambda}(z)}$ (Step 7). As
functions with summable coefficients are stable under
multiplication, this enables us to conclude.
\end{proof}

\emph{Step 9: $(z-1)(R(z)-I)^{-1}$ can be continuously
  extended on all $\overline{\D}$, and its Fourier coefficients (on
  $S^1$) are summable.}

\begin{proof}
We have already proved that $(z-1)(\tilde{R}(z)-I)^{-1}$ can be
continuously extended to $1$ on $S^1$. As $R$ and $\tilde{R}$
coincide in a neighborhood of $1$, it shows that
$(z-1)(R(z)-I)^{-1}$ can be continuously extended to $1$ on $S^1$.
Since we are interested in an extension to the whole disc
$\overline{\D}$, we must check that the previous arguments work
well on $\overline{\D}$, which does not present any difficulty:
dropping the tildes, Equation \eqref{exprime_valeur_propre} is
valid for $z$ in a neighborhood of $1$ in $\overline{\D}$, whence
$\frac{1-\lambda(z)}{1-z}$ tends to $\mu$ when $z\to 1$ in
$\overline{\D}$; using Equation \eqref{prolongement_en_1}, this
gives the desired extension to $1$.

On $S^1$,
    \begin{equation*}
    (z-1)(R(z)-I)^{-1}
    =(z-1)(\tilde{R}(z)-I)^{-1} \cdot
    (\tilde{R}(z)-I)(R(z)-I)^{-1}
    \end{equation*}
Step 8 shows that $(z-1)(\tilde{R}(z)-I)^{-1}$ has its Fourier
coefficients summable. Moreover, $(\tilde{R}(z)-I)(R(z)-I)^{-1}$
is $C^{1+\alpha}$ outside of any neighborhood of $1$, and equal to
$I$ on a neighborhood of $1$. Hence, it is $C^{1+\alpha}$ on
$S^1$ and has its coefficients summable. To conclude, we apply
Proposition \ref{stable_produit} which tells that the product of
functions with summable Fourier coefficients still has summable
coefficients.
\end{proof}

\emph{Step 10: $\sum \norm{T_{n+1}-T_n}<+\infty$.}
\begin{proof}
Let $A(z)=(1-z)(I-R(z))^{-1}$. For $|z|<1$, $A(z)=\sum
(T_n-T_{n-1})z^n$, so, when $r<1$, $T_n-T_{n-1}=\frac{1}{2\pi
r^n}\int_{0}^{2\pi} A(r e^{i\theta})e^{-i n \theta} \dd\theta$. As
$A$ can be continuously extended on $\overline{\D}$, we can let
$r$ tend to $1$ and obtain $T_n-T_{n-1}=c_n(A)$. But we have
already proved in the previous step that the coefficients of $A$
were summable.
\end{proof}

\section{Proof of the main theorem}
\label{section_preuve_main}
Once we have obtained the first main
lemma, the rest of the proof of Theorem
\ref{theoreme_renouvellement} is very similar to Sarig's
arguments. We will reproduce here only the parts which need to be
modified to fit in the current context.

To obtain the asymptotic expansion of $T_n$, the main idea is to
write $T(z)=\frac{1}{1-z}S(z)^{-1}$, where
$S(z)=\frac{I-R(z)}{1-z}$, to decompose $S=S_B+(S-S_B)$ where
$S_B(z)$ is a well controlled polynomial and $S-S_B$ a small
remainder, and to make a perturbative development of $S^{-1}$
using this decomposition. This amounts to writing
    \begin{equation}
    \label{devpt_perturbatif}
    T(z)=\frac{1}{1-z}S_B^{-1}+\frac{1}{1-z}S_B^{-1}(S_B-S)S_B^{-1}
    +\frac{1}{1-z}[S_B^{-1}(S_B-S)]^2 S^{-1}
    \end{equation}
The term $\frac{1}{1-z}S_B^{-1}(z)$ will give the contribution
$\frac{1}{\mu}P$ in the expansion of $T_n$, while the second one
will give the term
$\frac{1}{\mu²}\sum_{k=n+1}^{+\infty}P_k$ and the
third one will give the error term.

Write $S_B$ as $\frac{I-R_B(z)}{1-z}$ where
$R_B(z)=\sum_{n=1}^N z^n R_n+\sum_{n=N+1}^\infty R_n +(z-1)
\sum_{n=N+1}^\infty n R_n$ : this expression is such that
$R_B(1)=R(1)$ and $R'_B(1)=R'(1)$. For Equation
\eqref{devpt_perturbatif} to be valid for $z\in
\overline{\D}-\{1\}$, we have to check that $S_B$ is invertible,
i.e. $I-R_B$ is invertible. Following \cite[Proof of the Second
Main Lemma]{sarig:decay}, this is implied by the first main lemma
proved in the previous section as soon $N$ is large enough.

We recall without proof Sarig's second main lemma, which is a
consequence of the first main lemma.
\begin{lem}[Second Main Lemma]
\label{second_main_lemma}
Under the assumptions of Theorem
\ref{theoreme_renouvellement}, if $P$ is the eigenprojection of
$R(1)$ at $1$ and $\mu$ is given by $PR'(1)P=\mu P$, then there
exists $R_B:\C\to \Hom(\boL,\boL)$ with the following
properties:
\begin{enumerate}
\item $R_B$ is holomorphic, $R_B(1)=R(1)$ and $R'_B(1)=R'(1)$.
\item $\frac{R(1)-R_B}{1-z}$ and $\frac{1}{1-z}\left[ \frac{R(1)-
R_B}{1-z} -R'(1)\right]$ are polynomials in $z$.
\item $I-R_B(z)$ has a bounded inverse in $\Hom(\boL,\boL)$ for
every $z\in \overline{\D}-\{1\}$.
\item $\forall z\in\D$, $\left(\frac{I-R_B}{1-z}\right)^{-1}
= \frac{1}{\mu}P +(1-z)\sum_{n\geq 0} z^n A_n$ where $\norm{A_n}=
O( \kappa^n)$ for some $0<\kappa<1$.
\end{enumerate}
\end{lem}

Equation \eqref{devpt_perturbatif} together with the following
lemma (extending Sarig's Lemma 7 to the case $1<\beta\leq 2$ and
sharpening it for $\beta>2$) gives Theorem
\ref{theoreme_renouvellement}.
\begin{lem}
\label{estimees_sur_Tn}
Under the assumptions of Theorem
\ref{theoreme_renouvellement}, if $P$ is the eigenprojection of
$R(1)$ at $1$ and $\mu$ is given by $PR'(1)P=\mu P$, then
\begin{enumerate}
\item $\frac{1}{1-z}S_B^{-1}=\frac{1}{\mu}\sum_{n\geq 0} z^n(P+
\epsilon_n)$ where $\norm{\epsilon_n}=O(\kappa^n)$ for some
$0<\kappa<1$.
\item $\frac{1}{1-z}S_B^{-1}(S_B-S)S_B^{-1}=\frac{1}{\mu^2}
\sum_{n\geq 0}z^n \bigl(\sum_{k>n} P_k+\epsilon'_n\bigr)$ where
$\norm{\epsilon'_n}=O(1/n^\beta)$ and $P_n=\sum_{l>n}P R_l P$.
\item $\frac{1}{1-z}[S_B^{-1}(S_B-S)]^2 S^{-1}=\sum_{n\geq 0}z^n
E_n$ where
    $\norm{E_n}=\left\{\begin{array}{ll}
    O(1/n^\beta)&\text{if }\beta>2\\
    O(\log n/n^2)&\text{if }\beta=2\\
    O(1/n^{2\beta-2})&\text{if }2>\beta>1
    \end{array}\right.$
\end{enumerate}
\end{lem}

To prove the estimates in Lemma \ref{estimees_sur_Tn},
we will need some results on the
convolution of sequences. If $a_n$ and $b_n$ are sequences, put
$c_n=\sum_{k+l=n}a_kb_l$. We write $c=a\star b$.
\begin{lem}
\label{convole_suites} If $a_n=O(1/n^\alpha)$ and
$b_n=O(1/n^\beta)$ for some $\alpha\leq \beta\in\R$, then
    \begin{equation}
    (a\star b)_n=
    \left\{\begin{array}{ll}
    O(1/n^\alpha)&\text{if }\beta>1\\
    O(\log n/n^\alpha) &\text{if }\beta=1\\
    O(1/n^{\alpha+\beta-1}) & \text{if }\beta<1
    \end{array}\right.
    \end{equation}
In particular, for $\alpha>1$ or $\beta>1$ (without assuming
$\alpha\leq \beta$), $(a \star b)_n=O(1/n^\alpha)+O(1/n^\beta)$.
\end{lem}
\begin{proof}
We prove the result for $\beta<1$, the other cases being treated
in the same way. If $c_n=\sum_{k=0}^n a_k b_{n-k}$, we have
    \begin{equation*}
    |c_n|\leq \Bigl(\max_{0\leq k<n/2}|b_{n-k}|\Bigr) \sum_{0\leq k<n/2}|a_k|
    + \Bigl(\max_{n/2\leq k\leq n}|a_k|\Bigr) \sum_{n/2\leq k\leq n}|b_{n-k}|
    \end{equation*}
The sums can be estimated from above by $O(1/n^{\alpha-1})$ and
$O(1/n^{\beta-1})$ respectively, while the maxima are
$O(1/n^\beta)$ and $O(1/n^\alpha)$. This gives the conclusion.
\end{proof}
Let us state another lemma which will be useful later in
Section \ref{section_termes_suivants}. Its proof,
based on the same ideas, will be omitted.
\begin{lem}
\label{convole_suites_log}
If $a_n=O(\frac{(\log n)^u}{n^\alpha})$
and $b_n=O(\frac{(\log n)^v}{n^\beta})$ for some $\alpha\leq
1,\beta \leq 1$ and $u,v\geq 0$, then $(a\star b)_n=O(\frac{(\log
n)^{u+v+1}}{n^{\alpha+\beta-1}})$.
\end{lem}
In fact, the $(\log n)^{u+v+1}$ can be replaced by $(\log n)^{u+v}$
whenever $\alpha<1$ and $\beta<1$, but we will not need it.

We recall a notation used by Sarig: if $c_n$ is a real sequence
and $F(z)=\sum F_n z^n$ a formal series with coefficients in a
Banach algebra, write $F\in \Re(c_n)$ if $\norm{F_n}=O(c_n)$.
Abusing slightly notation, we write $\Re(1/n^\alpha)$ instead of
$\Re(1/(n+1)^\alpha)$, discarding the problem for $n=0$.

To prove Lemma \ref{estimees_sur_Tn}, we will first show that
$S(z)^{-1}\in \Re(1/n^\beta)$. In his main theorem, Sarig obtains
$\lfloor \beta \rfloor$ instead of $\beta$ since he proves only
that $S^{-1}\in \Re(1/n^{\lfloor \beta\rfloor})$; we can avoid
this loss of information with the help of Lemma
\ref{abstract_product} which should replace the general result on
Banach spaces Sarig uses and will give indeed $S(z)^{-1}\in
\Re(1/n^\beta)$.

\begin{lem}
\label{abstract_product}
Let $\boB$ be a Banach algebra and suppose that $F(z)=\sum F_k z^k$
where $\norm{F_k}=O(1/n^\beta)$ for some $\beta>1$. Suppose further that for
every $z\in\D$, $I+F(z)$ is invertible, and that $(I+F(z))^{-1}=\sum z^k
G_k$. If $\sum \norm{G_k}<\infty$, then $\norm{G_k}=O(1/n^\beta)$.
\end{lem}

Let us explain how to derive $S(z)^{-1}\in \Re(1/n^\beta)$ from this lemma.
Following Sarig, we use the identity
$S^{-1}=S_B^{-1}(I+(S-S_B)S_B^{-1})^{-1}$. In order to get the
result for $S^{-1}$, it is enough to prove that
$(I+(S-S_B)S_B^{-1})^{-1} \in \Re(1/n^\beta)$ since we already
know that $S_B^{-1}\in \Re(\kappa^n)$ for some $\kappa<1$
(Lemma \ref{second_main_lemma}). Note that
$(I+(S-S_B)S_B^{-1})^{-1}=S_BS^{-1}=I+(S_B-S)S^{-1}$ has its
coefficients summable because this is the case for $S^{-1}$ (Lemma
\ref{first_main_lemma}) and for $S_B-S$ (because
$S_B-S=\frac{R(1)-R_B}{1-z}-\frac{R(1)-R}{1-z}$, the first term
being a polynomial and the second one in $\Re(1/n^\beta)$).
Moreover, Lemma \ref{second_main_lemma} gives that
$I+(S-S_B)S_B^{-1}\in \Re(1/n^\beta)$ (since $S-S_B\in \Re(1/n^\beta)$ and
$S_B^{-1}\in \Re(\kappa^n)$). Consequently,
Lemma \ref{abstract_product} applied to $F=(S-S_B)S_B^{-1}$ yields
$(I+(S-S_B)S_B^{-1})^{-1} \in \Re(1/n^\beta)$, which gives
the conclusion.

\begin{proof}[Proof of Lemma \ref{abstract_product}]
Set $c_n=\sum_{i+k=n}\norm{G_i}\norm{G_k}$. As $\norm{G_n}$ is
summable, this is also the case for $c_n$. We will write $f_n$ and
$g_n$ respectively for $\norm{F_n}$ and $\norm{G_n}$.

Equating coefficients in
$\bigl[(I+F)^{-1}\bigr]'=-(I+F)^{-1}F'(I+F)^{-1}$ gives
    \begin{equation*}
    ng_n\leq \sum_{i+j+k=n}g_ijf_jg_k=\sum_{j=0}^n jf_jc_{n-j}
    \leq (\sup jf_j)\sum c_k<+\infty
    \end{equation*}
Consequently, $g_n=O(1/n)$. Moreover, we have $(ng_n)\leq c \star (jf_j)$, with
$jf_j=O(1/n^{\beta-1})$.

We show that $g_n=O(1/n^{1+\delta})$ for some $\delta>0$. It is
enough to prove this when $1<\beta<2$. As $g_n=O(1/n)$, $c=g\star
g$ is such that $c_n=O(\log n/n)$ according to Lemma
\ref{convole_suites}. Hence, $c_n=O(1/n^\gamma)$ for every $\gamma
<1$. Lemma \ref{convole_suites} again gives $c\star
(jf_j)=O(1/n^{\gamma+\beta-1-1})$, and
$g_n=O(1/n^{\gamma+\beta-1})$. As $\beta-1>0$ and $\gamma$ can be
chosen arbitrarily close to $1$, we can impose
$\gamma+\beta-1>1$, which gives the conclusion.

Assume that $g_n=O(1/n^\eta)$ for some $\eta>1$. As $c=g\star g$,
we get $c_n=O(1/n^\eta)$, whence $c\star (jf_j)=O(1/n^\eta)+
O(1/n^{\beta-1})$ once again by Lemma \ref{convole_suites}. As
$(ng_n)\leq c\star (jf_j)$, this implies $g_n=O(1/n^{\eta+1})+
O(1/n^\beta)$.

We already know that $g_n=O(1/n^{1+\delta})$ for some $\delta>0$.
Using the previous paragraph, we show by induction that, for any
integer $k$ such that $1+\delta+k<\beta$, we have
$g_n=O(1/n^{1+\delta+k+1})+ O(1/n^\beta)$. For the largest $k$
such that $1+\delta+k<\beta$, we obtain $g_n=O(1/n^\beta)$.
\end{proof}

From this point on, we can strictly follow Sarig's proof,
replacing his estimates $O(1/n^{\lfloor \beta \rfloor})$ by
$O(1/n^\beta)$. This way, we can obtain Estimates $(1)$ and $(2)$ in
Lemma \ref{estimees_sur_Tn}. However, the proof of Estimate $(3)$
has to be adapted.

\begin{proof}[Proof of Estimate $(3)$ in Lemma
\ref{estimees_sur_Tn}]
\label{preuve_estimee_3}
As in \cite[Step 4 of the proof of Theorem 1]{sarig:decay}, write
$G(z)=S_B^{-1}(z)(S_B(z)-S(z))=\sum G_k z^k$. As $S_B-S\in \Re(1/n^\beta)$
and $S_B^{-1}\in \Re(\kappa^n)$ for some $\kappa<1$
(Lemma \ref{second_main_lemma}), we obtain that $G\in \Re(1/n^\beta)$.
Moreover, $\sum G_k=0$ (because $S_B(1)=S(1)$), hence
$\frac{1}{1-z}G(z)=-\sum z^n\sum_{k>n}G_k$ and consequently
$\frac{G(z)}{1-z}\in \Re(1/n^{\beta-1})$ (see \cite[Step 4 of the
proof of Theorem 1]{sarig:decay} for more details).

Setting $E=\frac{1}{1-z}G^2 S^{-1}=\sum z^n E_n$, we want to
estimate the coefficients $E_n$ of $E$. We have
    \begin{equation*}
    E'=\left(\frac{G}{1-z}\right)^2
    S^{-1}+\left[\left(\frac{G}{1-z}\right)G'+G'
    \left(\frac{G}{1-z}\right)\right]S^{-1}+G\left(\frac{G}{1-z}\right)
    (S^{-1})'
    \end{equation*}
We know that $\frac{G}{1-z}\in \Re(\frac{1}{n^{\beta-1}})$, $G'\in
\Re(\frac{1}{n^{\beta-1}})$ and $(S^{-1})'\in
\Re(\frac{1}{n^{\beta-1}})$ (since $S^{-1}\in
\Re(\frac{1}{n^\beta})$). Lemma \ref{convole_suites} on
convolutions gives
    \begin{equation*}
    \left(\frac{G}{1-z}\right)^2\in \left\{\begin{array}{ll}
    \Re(1/n^{\beta-1})&\text{if }\beta>2\\
    \Re(\log n/n) &\text{if }\beta=2\\
    \Re(1/n^{2\beta-3}) &\text{if }\beta<2
    \end{array}\right.
    \end{equation*}
and we have analogous estimates for the other terms in $E'$.
Integrating, we get the desired estimates for $E_n$. This
concludes the proof of Lemma \ref{estimees_sur_Tn}
and, with it, of Theorem \ref{theoreme_renouvellement}.
\end{proof}

\section{Higher order terms in $T_n$\label{section_termes_suivants}}

To obtain an asymptotic expansion of $T_n$, we have used the
perturbative development of $T(z)$ up to order 2 given in Equation
\eqref{devpt_perturbatif}.
While this is enough to obtain the asymptotic expansion of $T_n$
stated in Theorem \ref{theoreme_renouvellement}, this does not
give an optimal result: in the course of the proof, we have shown
that $T_{n+1}-T_n=O(1/n^\beta)$, whence an expansion up to order
$1/n^\beta$ can be expected. In this section, we will show that it
is indeed possible to obtain this estimate and we will give the
first terms of this expansion.

Fix $N\in \N$ the order of the expansion. Then we have
    \begin{equation}
      \label{Developpement_perturbatif_Tn}
    T(z)=\frac{1}{1-z}\sum_{k=0}^{N-1}[S_B^{-1}(S_B-S)]^k S_B^{-1}+
    \frac{1}{1-z}[S_B^{-1}(S_B-S)]^N S^{-1}
    \end{equation}
To prove that this equation gives good asymptotics on $T_n$, we
have to see that the coefficients of
$\frac{1}{1-z}[S_B^{-1}(S_B-S)]^N S^{-1}$ decrease fast to zero,
at least at a speed $O(1/n^\beta)$, if $N$ is large enough. We
will use the fact that, if $G(z)=S_B^{-1}(S_B-S)$, then $G\in
\Re(1/n^\beta)$ and $G/(1-z)\in \Re(1/n^{\beta-1})$, as we have
seen in the proof of Estimate $(3)$ in Lemma
\ref{estimees_sur_Tn}.

\begin{lem}
\label{lemme_technique_induction}
Let $G(z)=\sum G_n z^n$ be a formal series with coefficients in a Banach
algebra $\boB$, such that $G(z)\in \Re(1/n^\beta)$ and $\frac{G(z)}{1-z}\in
\Re(1/n^{\beta-1})$ for some $\beta>1$. Then, for any $p\in \N$, there exists a
constant $C$
such that for any $H_1,\ldots,H_{p-1}\in \boB$,
    \begin{equation*}
    \norm{ \left(\frac{G(z)H_1G(z)\ldots H_{p-1}G(z)}{1-z}\right)_n}
    \leq C \norm{H_1}\ldots\norm{H_{p-1}} \cdot
    \left\{\begin{array}{ll}
        \frac{1}{n^\beta} &\text{if }\beta<p(\beta-1)\\
        \frac{\log n}{n^\beta} &\text{if }\beta=p(\beta-1)\\
        \frac{1}{n^{p(\beta-1)}} &\text{if }\beta>p(\beta-1)
    \end{array}\right.
    \end{equation*}
\end{lem}
(the notation $(\ \ )_n$ denotes the coefficient of $z^n$ in the
formal series between the braces).

By Lemma \ref{convole_suites}, the convolution of two sequences in
$\Re(1/n^\gamma)$ and $\Re(1/n^\delta)$ does not cause a serious
loss of information when $\gamma\leq 1$ and $\delta\leq 1$ (we
get a term in $\Re(1/n^{\gamma+\delta-1})$ with possibly a $\log
n$ if $\gamma$ or $\delta=1$), while there will be some loss of
information if one of the exponents is $>1$ (because a minimum
appears). The idea of the proof of Lemma
\ref{lemme_technique_induction} is to differentiate, which will
give exponents less than $1$ if $\beta < 2$. The problem is then
to control the terms which have not been differentiated -- this
is done using the $1/(1-z)$ and the induction.

\begin{proof}[Proof of Lemma \ref{lemme_technique_induction}]
We prove Lemma \ref{lemme_technique_induction} by induction on
$p$. The result is part of the hypotheses when $p=1$ and easily
follows from the proof of Estimate (3) in Lemma
\ref{estimees_sur_Tn} for $p=2$ (the same argument works when a
term $H_1$ is inserted). Assume $p>2$. For the moment, we will
forget about the $H_i$'s and see later that they do not matter.

If $\beta<(p-1)(\beta-1)$, the induction gives that
  $\frac{G(z)^{p-1}}{1-z}\in \Re(1/n^\beta)$. As $G(z)\in \Re(1/n^\beta)$, a
simple convolution gives the desired result. Thus, we can assume
that $\beta\geq(p-1)(\beta-1)$. As $p\geq 3$, this implies in
particular that $\beta \leq 2$.

Differentiating $p-1$ times $F(z)=\frac{G(z)^p}{1-z}$ gives, for
some constants $C_{i,i_1,\ldots,i_p}$
    \begin{equation*}
    F^{(p-1)}(z)=\sum_{\substack{i+i_1+\ldots+i_p=p-1\\ i_1,\ldots,i_p\geq 0}}C_{i,i_1,\ldots,i_p}\frac{G(z)^{(i_1)}\ldots
      G(z)^{(i_p)}}{(1-z)^{i+1}}
    \end{equation*}
where  $G(z)^{(k)}$ denotes the function $G(z)$ differentiated $k$ times.
We will do the proof assuming that the $G(z)^{(k)}$ commute, and see later what
should be modified in the general case.

Permute the $G(z)^{(i_j)}$ and group at the end the coefficients
where $i_j=0$, i.e. the factors which are not differentiated.
    \begin{equation*}
    F^{(p-1)}(z)=\sum_{k=0}^{p-1}\sum_{\substack{i+i_1+\ldots+i_k=p-1\\
        i_1,\ldots,i_k\not=0}} D_{i,i_1,\ldots,i_k} G^{(i_1)}\ldots G^{(i_k)}
    \frac{G(z)^{p-k}}{(1-z)^{i+1}}
    \end{equation*}
As the $i_j$'s are nonzero, $i+k\leq p-1$, i.e. $i+1\leq p-k$.
Consequently, we write the factor $\frac{G^{p-k}}{(1-z)^{i+1}}$ as
$\left(\frac{G(z)}{1-z}\right)^i
  \frac{G(z)^{p-k-i}}{1-z}$. If $\beta<2$, the terms $G^{(i_j)}$ and $G/(1-z)$ are in
  $\Re(1/n^\gamma)$ for exponents $\gamma< 1$, which means that we will
not loose information when multiplying them. In fact, whenever
$\beta<2$, the product $G^{(i_1)}\ldots G^{(i_k)}
\left(\frac{G(z)}{1-z}\right)^i$ will be in $\Re(1/n^\delta)$
where $\delta=\sum_{t=1}^k (\beta-i_t) + i(\beta-1)-k-i+1$, using
 Lemma \ref{convole_suites} $k+i-1$ times with exponents $<1$, and
if $\beta=2$ it will be in $\Re(\frac{(\log n)^u}{n^\delta})$ for
some integer $u$, by Lemma \ref{convole_suites_log}. Simplifying,
$\delta=(k+i)(\beta-1)-p+2$. The remaining factor
$\frac{G(z)^{p-k-i}}{1-z}$ will be controlled using the induction.

Let us distinguish 3 cases:
\begin{enumerate}
\item{If $(p-1)(\beta-1)\leq \beta < p(\beta-1)$.}

For every term $S_{i,i_1,\ldots,i_k}=G^{(i_1)} \ldots G^{(i_k)}
\left(\frac{G}{1-z}\right)^i \frac{G^{p-k-i}}{1-z}$, we have
$p-k-i<p$. Consequently, $(p-k-i)(\beta-1)\leq \beta$, and we are
in the second or third case of the induction. In fact, we are in
the second case only if $(p-k-i)(\beta-1)= \beta$, which implies
that $k=1$ and $i=0$.

Let us first consider this term corresponding to $i=0, i_1=p-1$:
$S_{0,p-1}=G(z)^{(p-1)} \frac{G(z)^{p-1}}{1-z}$. The induction
hypothesis ensures that $\frac{G(z)^{p-1}}{1-z}\in \Re(\frac{\log
  n}{n^{(p-1)(\beta-1)}})$. As
$p(\beta-1)>\beta$, we have $(p-1)(\beta-1)>1$. Thus, $\frac{\log
  n}{n^{(p-1)(\beta-1)}}=O(1/n^\gamma)$ for some $\gamma>1$, and
$\frac{G(z)^{p-1}}{1-z}\in \Re(\frac{1}{n^\gamma})$.  When convolving with
$G(z)^{(p-1)}\in \Re(\frac{1}{n^{\beta-p+1}})$ (where
$\beta-p+1\leq 1$ since $\beta\leq 2$), Lemma
\ref{convole_suites} gives an expression in
$\Re(\frac{1}{n^{\beta-p+1}})$.

Now, we consider the other terms, of the form
$S_{i,i_1,\ldots,i_k}=G^{(i_1)} \ldots G^{(i_k)}
\left(\frac{G}{1-z}\right)^i \frac{G^{p-k-i}}{1-z}$ with $k+i>1$.
As $(p-1)(\beta-1)\leq \beta$, we obtain $(p-2)(\beta-1)\leq 1$,
hence $(p-k-i)(\beta-1)\leq 1$. As $\frac{G^{p-k-i}}{1-z} \in
\Re(\frac{1}{n^{(p-k-i)(\beta-1)}})$ by induction and $G^{(i_1)}
\ldots G^{(i_k)} \left(\frac{G}{1-z}\right)^i \in \Re(\frac{(\log
n)^u }{n^\delta})$ where $\delta=(k+i)(\beta-1)-p+2$, we can
convolve and obtain $S_{i,i_1,\ldots,i_k}\in \Re(\frac{(\log
n)^{u+1}}{n^{\nu}})$ where
$\nu=\delta+(p-k-i)(\beta-1)-1=p(\beta-1)-(p-1)$. Since
$p(\beta-1)>\beta$, we have in fact $\nu>\beta-(p-1)$, which
implies that $S_{i,i_1,\ldots,i_k}\in
\Re(\frac{1}{n^{\beta-(p-1)}})$.

Summing all terms, we obtain $F(z)^{(p-1)}\in
\Re(\frac{1}{n^{\beta-p+1}})$. Integrating then $p-1$ times, we
get $F(z)\in \Re(\frac{1}{n^\beta})$, which is the desired result.

\item{If $\beta=p(\beta-1)$.}

Here, we have $\beta<2$, which implies that the term $(\log n)^u$
disappears.

We use the same reasoning as in case 1. Here, $(p-1)(\beta-1)=1$,
which means that, to obtain $S_{0.p-1}$, we have to convolve
elements in $\Re(\frac{1}{n^{\beta-p+1}})$ and in
$\Re(\frac{1}{n})$, which gives $\Re(\frac{\log
  n}{n^{\beta-p+1}})$. The other terms $S_{i,i_1,\ldots,i_k}$ are treated as
above and are in $\Re(\frac{\log n}{n^{p(\beta-1)-p+1}})=\Re(\frac{\log
  n}{n^{\beta -p+1}})$. Summing and integrating gives the result.

\item{If $\beta>p(\beta-1)$.}

We do not need to distinguish the term $S_{0,p-1}$ any more: in all
terms, all exponents are $<1$. A convolution gives terms in
$\Re(\frac{1}{n^{p(\beta-1)-p+1}})$, which gives the result after integration.
\end{enumerate}

This concludes the proof in the commutative case, and without the
$H_j$. In fact, in the commutative case, there would be no more
trouble to include the $H_j$, since we can group them for example
at the beginning and write $\frac{G(z)H_1G(z)\ldots
H_{p-1}G(z)}{1-z}=H_1\ldots H_{p-1} \frac{G(z)^{p-1}}{1-z}$; the
result proved on $\frac{G(z)^{p-1}}{1-z}$ enables us to conclude.

What remains to be done is to see how to get rid of this
commutativity hypothesis for $p\geq 3$. To avoid cumbersome
notation, we will only see on an example what happens.

Suppose that we want to estimate $F(z)=\frac{G(z)^3}{1-z}$. After
two differentiations, we obtain many terms that can be controlled
using the previous method (since the factors $\frac{G(z)^d}{1-z}$
are already grouped), and one term $\frac{G(z)G''(z) G(z)}{1-z}$.
This term is problematic: what we would like to do is to convolve
first the two extremal $G(z)$ and $1/(1-z)$, to use the induction
hypothesis to see that this is of the order of
$\frac{G(z)^2}{1-z}$, i.e. of order $O(1/n^{2\beta-2})$, and then
to convolve it with $G''(z)$ to get the result. This is indeed
possible since, if we write $G''(z)= \sum J_n z^n$, we have
    \begin{equation*}
    \left(\frac{G(z)G''(z)G(z)}{1-z}\right)_n=\sum_{k+l=n} \left(\frac{G(z)J_l G(z)}{1-z}\right)_k
    \end{equation*}
The induction hypothesis (with the $H_l$ replaced by the $J_l$)
gives a control in $O(1/n^{2\beta-2})$ on each term of the
previous sum. We obtain
    \begin{equation*}
    \norm{\left(\frac{G(z)G''(z)G(z)}{1-z}\right)_n}
    \leq \sum_{k+l=n} C\norm{J_l} \frac{1}{k^{2\beta-2}}
    \end{equation*}
which amounts to the usual convolution between $\norm{J_l}$ and
$1/n^{2\beta-2}$. This shows that, using the induction
hypothesis, we are able to obtain the same estimate on
$\frac{G(z)G''(z)G(z)}{1-z}$ as on $\frac{G(z)^2 G''(z)}{1-z}$. It
can indeed be done for as many $G(z)$ terms as necessary, which
shows that all the previous estimates in the commutative case
apply also in the general case.
\end{proof}

\begin{lem}[Control of the error term]
\label{controle_error_term}
We have
    \begin{equation*}
    \frac{1}{1-z}[S_B^{-1}(S_B-S)]^N S^{-1} \in \left\{
    \begin{array}{cc}
    \Re\left(\frac{1}{n^\beta}\right) &\text{if }N(\beta-1)>\beta\\
    \Re\left(\frac{\log n}{n^\beta}\right) &\text{if }N(\beta-1)=\beta\\
    \Re\left(\frac{1}{n^{N(\beta-1)}}\right) &\text{if }N(\beta-1)<\beta
    \end{array}\right.
    \end{equation*}
\end{lem}
\begin{proof}
Set $G(z)=S_B^{-1}(S_B-S)$. Then the conditions of
Lemma \ref{lemme_technique_induction} are verified (this has
been checked in the proof of
Estimate (3) in Lemma \ref{estimees_sur_Tn}).
Consequently, the lemma gives
estimates on $\frac{G(z)^N}{1-z}$. As we already know that
$S(z)^{-1}\in \Re(1/n^\beta)$ with $\beta>1$, another convolution
enables us to conclude.
\end{proof}

To use this result, there remains only to study the terms in the
sum in Equation \eqref{Developpement_perturbatif_Tn}, i.e. the
perturbative terms for $k=0,\ldots,N-1$. The method used in Sarig
to estimate the first term still works: estimating $S_B^{-1}$ by
$\frac{1}{\mu}P$ gives an exponentially decreasing error, which
does not matter. Moreover, we can estimate $S_B-S$ by
$\sum_{k=0}^{\infty} (1-z^k)\sum_{n=k+1}^\infty R_k$. A formal
multiplication gives finally the desired terms. More precisely,
the following lemma is valid.

\begin{lem}[Estimates on the perturbative terms]
\label{lemme_estimee_perturb}
\ \\
For any $k\in \N^*$, writing $P_n=\sum_{l>n}PR_lP$, we have
    \begin{equation*}
    \frac{1}{1-z}[S_B^{-1}(S_B-S)]^k S_B^{-1}=\frac{1}{\mu^{k+1}} \frac{1}{1-z}\left[
    \sum_{n=0}^\infty (1-z^n)P_n\right]^k+E(z)
    \end{equation*}
where $E\in \Re(1/n^\beta)$.
\end{lem}
\begin{proof}
We already know that $S_B-S \in \Re(1/n^\beta)$ and $S_B^{-1}$
also.

We write
    \begin{align*}
    S_B-S=\left[R'(1)-\frac{R(1)-R}{1-z}\right]+
    \left[\frac{R(1)-R_B}{1-z}-R'(1)\right]
    \\
    =\left[\sum_{n=0}^\infty (1-z^n) \sum_{l=n+1}^\infty R_l
    \right]
    +(1-z)B(z)
    \end{align*}
where $B(z)$ is a polynomial, according to Lemma
\ref{second_main_lemma}. Moreover, the same lemma gives that
    \begin{equation*}
    S_B^{-1}=\frac{1}{\mu}P+(1-z)A(z)
    \end{equation*}
for some $A(z)\in \Re(\kappa^n)$ with $\kappa <1$.

We multiply these expressions to get
$\frac{1}{1-z}[S_B^{-1}(S_B-S)]^k S_B^{-1}$ and we expand the
product. If we choose a term $(1-z)A(z)$ or $(1-z)B(z)$, we use
it to simplify the $\frac{1}{1-z}$, and all the other terms are
$\Re(1/n^\beta)$, which gives after convolution still a
$\Re(1/n^\beta)$. The remaining term gives the expression stated
in the lemma.
\end{proof}

Gathering the results of Lemma \ref{controle_error_term} and
Lemma \ref{lemme_estimee_perturb}, we obtain
\begin{thm}
\label{Tn_ordre_plus_grand} Under the hypotheses of Theorem
\ref{theoreme_renouvellement}, we have, for any $N\in \N$,
writing $P_m=\sum_{k>m}PR_kP$,
    \begin{equation}
    \label{controle_erreur}
    \sum T_n z^n=\frac{1}{1-z}\frac{1}{\mu}P+\sum_{k=1}^{N-1}\frac{1}{\mu^{k+1}} \frac{1}{1-z}\left[
    \sum_{m=0}^\infty (1-z^m)P_m\right]^k
    +E(z)
    \end{equation}
where
    \begin{equation*}
    E(z)\in \left\{
    \begin{array}{cc}
    \Re\left(\frac{1}{n^\beta}\right) &\text{if }N(\beta-1)>\beta\\
    \Re\left(\frac{\log n}{n^\beta}\right) &\text{if }N(\beta-1)=\beta\\
    \Re\left(\frac{1}{n^{N(\beta-1)}}\right) &\text{if }N(\beta-1)<\beta
    \end{array}\right.
    \end{equation*}
\end{thm}
Note that, for any $\beta>1$, it is possible to choose $N$ such
that $N(\beta-1)>\beta$, which implies that the expansion of $T_n$
with $N$ terms gives an estimate with an error term in
$O(1/n^\beta)$. In particular, if $Pf=0$, we obtain $T_n
f=O(1/n^\beta)$, which is exactly Theorem \ref{thm_quand_Pf_nul}.

To obtain a sharp asymptotic expansion for $T_n$, there remains
only to expand the middle terms in Equation
\eqref{controle_erreur}. We give for example the theorem that we
obtain for $N=3$:

\begin{thm}
Under the hypotheses of Theorem \ref{theoreme_renouvellement}, we
have
    \begin{equation*}
    T_n=\frac{1}{\mu}P+\frac{1}{\mu^2}\sum_{k=n+1}^{+\infty}P_k
    +\frac{1}{\mu^3}\Bigl( \sum_{k,l>n}P_k P_l - \sum_{\substack{k,l\leq
    n\\k+l>n}} P_kP_l\Bigr)+E_n
    \end{equation*}
where $E_n\in \Hom(\boL,\boL)$ satisfy
    \begin{equation*}
    \norm{E_n}=\left\{\begin{array}{ll}
    O(1/n^\beta)&\text{if }\beta>3/2\\
    O(\log n/n^\beta)&\text{if }\beta=3/2\\
    O(1/n^{3(\beta-1)})&\text{if }3/2>\beta>1
    \end{array}\right.
    \end{equation*}
\end{thm}

We give for completeness the next term in the expansion: after
tedious calculations, we find that it is (up to the factor
$1/\mu^4$)
    \begin{gather*}
    \sum_{k,l,m>n}
    -\sum_{\substack{0<k,l\leq n\\k+l>n\\m>n}}
    -\sum_{\substack{0<k,m\leq n\\k+m>n\\l>n}}
    -\sum_{\substack{0<l,m\leq n\\l+m>n\\k>n}}
    -\sum_{\substack{0<k,l,m\leq n\\k+l>n\\k+m>n}}
    \\
    -\sum_{\substack{0<k,l,m\leq n\\l+k>n\\l+m>n}}
    -\sum_{\substack{0<k,l,m\leq n\\m+k>n\\m+l>n}}
    +\sum_{\substack{0<k,l,m\leq n\\k+l>n\\k+m>n\\l+m>n}}
    +\sum_{\substack{k+l\leq n\\k+m\leq n\\l+m\leq n\\k+l+m>n}}
    P_kP_lP_m
    \end{gather*}

\section{Application to Markov maps \label{markov_maps}}

\subsection{Definition of Markov maps}
The definitions and results of this section are for the main part contained in
\cite{aaronson:book}.

A Markov map is a non-singular transformation $T$ of a Lebesgue
space $(X,\boB,m)$ together with a measurable partition $\alpha$
of $X$ such that if $a\in \alpha$, $m(a)>0$, $Ta$ is a union (mod
$m$) of elements of $\alpha$, and $T:a\to Ta$ is invertible.
Moreover, it is assumed that $\bigvee_0^\infty
T^{-i}\alpha=\boB$, i.e. the partition separates the points.

For $a_0,\ldots,a_{n-1}\in \alpha$ define a \emph{cylinder} by
$[a_0,\ldots,a_{n-1}]=\bigcap_{i=0}^{n-1}T^{-i}a_i$: two points
in a same cylinder of length $n$ remain in the same elements of
the partition up to time $n$. These cylinders can be used to
topologize the space $X$ (mod $0$), and even to define a metric
on it: $d_\theta(x,y)=\theta^{t(x,y)}$, where $t(x,y)=\sup\{n\tq
x,y\in
  [a_0,\ldots,a_{n-1}] \text{ for some }a_0,\ldots,a_{n-1}\in
  \alpha\}$ is the time until which $x$ and $y$ remain in the same
  elements of the partition $\alpha$, and $0<\theta<1$ is some fixed number.

A Markov map $T$ is said to be \emph{irreducible} if $\forall a,b\in
\alpha$, $\exists n$, $T^{-n}a \cap b\not=\emptyset$ (i.e. $b\subset
T^n a$). This means that there is no $\gamma \varsubsetneq \alpha$
such that the elements of $\gamma$ are stable by $T$.
An irreducible Markov map $T$ is \emph{aperiodic} if $\forall a\in
\alpha, \exists N\in \N, \forall n\geq N, a\subset T^n
a$. Equivalently, there exists such an $a$, or there exists an $a$
such that $\gcd\{n \tq a\subset T^n a\}=1$.
An irreducible aperiodic Markov map is also said to be
\emph{topologically mixing}, i.e. $\forall a,b \in \alpha, \exists N,
\forall n \geq N, b\subset T^n a$. This corresponds to the topological
mixing for the topology defined by the cylinders as above.

The transfer operator $\hat{T}$ associated to $T$ can be written
$\hat{T}f(x)=\sum_{Ty=x}g_m(y)f(y)$, where the weight $g_m$ is defined
by $g_m=\frac{\de m}{\de m\circ T}$. Different regularity assumptions
are possible on $\log g_m$, corresponding to different controls of the
distortion.

For any function $\phi:X \to \C$, the \emph{variations} of $\phi$ are
defined by $v_n(\phi)=\sup\{ |\phi(x)-\phi(y)| \tq x,y\in
[a_0,\ldots,a_{n-1}] \text{ where }a_i\in \alpha\} $. The function
$\phi$ is said to have summable variations if $\sum_{n\geq
  1}v_n(\phi)<+\infty$, and to be Hölder continuous for the exponent
$\theta$ if $\exists C>0, \forall n\geq 1, v_n(\phi)\leq C
\theta^n$ (this is a definition, which corresponds to being
Lipschitzian with respect to the metric $d_\theta$ on each
element of the partition $\alpha$ ).

If $\log g_m$ is of summable variations, the distortion is bounded,
meaning that there exists a constant $C$ such that, for all $x,y\in
[a_0,\ldots,a_{n-1}]$,
$\left|\frac{g_m^{(n)}(x)}{g_m^{(n)}(y)}-1\right| \leq C$, where
$g_m^{(n)}=\prod_{i=0}^{n-1} g_m \circ T^i$ is the weight associated to
$\hat{T}^n$. In particular, this implies that $g_m^{(n)}(x)=C^{±
  1}\frac{m[a_0, \ldots, a_{n-1}]}{m[Ta_{n-1}]}$ (i.e.
$\frac{1}{C} \frac{m[a_0, \ldots, a_{n-1}]}{m[Ta_{n-1}]}\leq
g_m^{(n)}(x)\leq C\frac{m[a_0, \ldots, a_{n-1}]}{m[Ta_{n-1}]}$).
When the \emph{``big
  image'' property} $\inf_{a\in \alpha} m[Ta]>0$ is satisfied, we
obtain even $g_m^{(n)}(x)=D^{± 1} m[a_0,\ldots,a_{n-1}]$.

When the distortion is of summable variations, the application $T$
behaves almost like a probabilistic Markov chain, with
independence. Hence, it is possible to prove theorems ensuring a
good behavior. In particular,

\begin{prop}
\label{hypotheses_donnent_exactitude} Let $(X,\boB,T,m,\alpha)$
be an irreducible Markov map with the big image property for
which $\log g_m$ is of summable variations. Then $T$ is
conservative and ergodic.
\end{prop}
\begin{proof}
This is a corollary of Theorem 4.6.3 in \cite{aaronson:book}
(where the hypotheses are in fact weaker, since this theorem
requires only the ``weak distortion property'').
\end{proof}

If $\log g_m$ is Hölder continuous, the distortion is better
controlled, which gives stronger results. In particular, the
transfer operator $\hat{T}$ acting on the space of Hölder
continuous bounded functions admits a spectral gap (\cite[Thm
4.7.7]{aaronson:book}). More precisely, let $\alpha'$ denote the
smallest partition such that, $\forall a\in \alpha, Ta$ is
$\alpha'$-measurable; the partition $\alpha'$ is coarser than
$\alpha$. For $a\in \alpha'$ and $f:X\to \R$, write $D_a
f=\sup\{|f(x)-f(y)|/d_\theta(x,y) \tq x,y\in a\}$ the best
Lipschitz constant of $f$ on $a$. Finally, let $\boL$ be the space
of functions $f:X\to \C$ such that
$\norm{f}_\boL=\norm{f}_\infty+\sup_{a\in \alpha'} D_a f
<+\infty$. It is the space of Lipschitzian functions on $X$, but
the norm is not the usual Lipschitz norm. When $\log g_m$ is
Hölder continuous (for some exponent $\theta$) and $T$ has the
big image property, Ruelle has proved that the essential spectral
radius of $T$ acting on $\boL$ verifies $r_{ess}(T)\leq \theta$.

\subsection{Induced Markov maps}
From this point on, $(X,\boB,m,T,\alpha)$ will be a probability preserving
Markov map.

Let $\emptyset\not=\gamma \subset \alpha$. If $Y=\bigcup \gamma$,
the induced map $T_\gamma:Y\to Y$ is defined as the first return
map from $Y$ to $Y$, i.e. $T_\gamma=T^{\phi_\gamma}$, where
$\phi_\gamma(x)=\inf\{ n\geq 1\tq T^n(x)\in Y\}$ is the return
time to $Y$. If $x\not\in Y$, we set $\phi_\gamma(x)=0$. By the
Poincaré recurrence theorem, $T_\gamma$ and all its iterates are
defined for $m$-almost every point of $Y$ -- replacing $Y$ by
this smaller set, we can assume that $T_\gamma$ is in fact
defined on all $Y$.

A measure $m_\gamma$ is defined on $Y$ by $m_\gamma=m_{|Y}$. As $m$ is
invariant par $T$, the measure $m_\gamma$ is invariant by $T_\gamma$.

Let $\delta=\{[a,\xi_1,\ldots,\xi_{n-1},\gamma] \tq a\in \gamma,
\xi_1,\ldots,\xi_{n-1}\not \in \gamma,
[a,\xi_1,\ldots,\xi_{n-1},\gamma]\not=\emptyset\}$: this is a
partition of $Y$, for which $T_\gamma$ is a Markov map. The
cylinders for this partition will be denoted by
$[d_0,\ldots,d_{n-1}]_\gamma$ (with $d_0,\ldots,d_{n-1}\in
\delta$). If $d=[a,\xi_1,\ldots,\xi_{n-1},\gamma]\in \delta$, its
image is $T_\gamma d=T \xi_{n-1}$ -- hence, it is
$\gamma$-measurable. In particular, if $\gamma$ is finite, its
elements have a measure $\geq \epsilon>0$, which implies that
$\forall d\in \delta$, $m_\gamma (T_\gamma d)\geq \epsilon$. Thus
$T_\gamma$ has the ``big image'' property.

The following straightforward lemma establishes a link between the
mixing properties of $T$ and those of the induced transformation $T_\gamma$.
\begin{lem}
\label{top_mixing_induit_bien}
If $T$ is irreducible, then
$T_\gamma$ is irreducible.
\end{lem}
%
%

We will be interested in induced maps which have good distortion
properties. More precisely, write $g_{m_\gamma}=\frac{\de m_\gamma}{\de m_\gamma
  \circ T_\gamma}$. We assume that there exist constants $C>0$ and $\theta<1$
such that $\forall n\geq 1, v_n(\log g_{m_\gamma}) \leq C\theta^n$ (where $v_n$ is the
variation with respect to the induced map $T_\gamma$): we say that $\log
g_{m_\gamma}$ is locally Hölder continuous. In this case, the previous
theorems on maps whose distortion has summable variations
apply to $T_\gamma$.

As above, let $\delta'$ denote the smallest partition such that
$\forall d\in \delta$, $T_\gamma d$ is a union of atoms of
$\delta'$. As every $T_\gamma d$ is $\gamma$-measurable, this
partition is coarser than $\gamma$. For $x,y\in Y$, let
$t_\gamma(x,y)=\sup\{ n \tq x,y\in [d_0,\ldots,d_{n-1}]_\gamma\}$
and let $\boL$ denote the space of functions $f:Y\to \C$ such that
$\norm{f}_\boL:= \norm{f}_\infty+\sup_{d\in \delta'} D_d f
<+\infty $, where $D_d f$ is the least Lipschitz constant of $f$
on $d$ for the distance $d(x,y)=\theta^{t_\gamma(x,y)}$.

We now state the main theorem of this section:
\begin{thm}
\label{thm_applications_markov_precis} Let $(X,\boB,m,T,\alpha)$
be a topologically mixing probability preserving Markov map, and
$\emptyset\not=\gamma \subset \alpha$. Assume that $T_\gamma$ has
the big image property and that $g_{m_\gamma}$ has a version such
that $\log g_{m_\gamma}$ is locally $\theta$-Hölder continuous
for some $0<\theta<1$. Assume moreover that
$m[\phi_\gamma>n]=O(1/n^\beta)$ for some $\beta>1$.

Then $\exists C>0$ such that $\forall f,g$ integrable and
supported inside $Y$,
    \begin{equation*}
    \left|\Cor(f,g\circ T^n)-\left(
    \sum_{k=n+1}^{\infty}m[\phi_\gamma>k]\right)\int f \int
    g \right| \leq C F_\beta(n) \norm{g}_\infty \norm{f}_{\boL}
    \end{equation*}
where $F_\beta(n)=1/n^\beta$ if $\beta>2$, $\log n/n^2$ if $\beta=2$ and
$1/n^{2\beta-2}$ if $2>\beta>1$ (and $\boL$ denotes the space of
$\theta$-Hölderian functions on $Y$).

Moreover, if $\int f=0$, then $\Cor(f,g\circ T^n)=O(1/n^\beta)$.
\end{thm}

\subsection{Proof of Theorem \ref{thm_applications_markov_precis}}
The strategy is to apply the abstract Theorem
\ref{theoreme_renouvellement} to ``first
return transfer operators''. In this section, $(T,\alpha)$ will be a Markov
map and $\gamma\subset \alpha$ a subpartition such that the hypotheses of Theorem
\ref{thm_applications_markov_precis} are verified. The first three lemmas can
essentially be found in \cite{sarig:decay}.

For $\underline{d}=[d_0,\ldots,d_{n-1}]_\gamma \not =\emptyset$, define
$M_{\underline{d}} f(x)=g_{m_\gamma}^{(n)}(x) f(d_0\ldots d_{n-1}x)$ if this
point is defined, $0$ otherwise.
\begin{lem}
\label{lemme_sur_Mdbarre}
There exists a constant $B$ such that, $\forall
\underline{d}=[d_0,\ldots,d_{n-1}]_\gamma$, $\forall f\in \boL$,
  \begin{equation*}
  \norm{M_{\underline{d}}f}_\boL \leq B m[\underline{d}]\left(\theta^n
    \norm{f}_\boL+ \frac{1}{m[\underline{d}]} \int_{[\underline{d}]} |f|
      \dd m \right)
  \end{equation*}
\end{lem}
\begin{proof}
This lemma is classical and uses the distortion control to obtain
explicit estimates. See for example \cite[Lemma 8]{sarig:decay} or
\cite[Lemma 4.7.2]{aaronson:book}.
\end{proof}

Let $L$ be the operator defined by $L
f(x)=\sum_{Ty=x}g_m(y)f(y)$: it is a version of the transfer
operator $\hat{T}$, but it acts on actual functions and not on
functions defined almost everywhere. In the same way, but for the
induced map, set $L_\gamma f(x)=\sum_{T_\gamma
  y=x} g_{m_\gamma}(y)f(y)$. Write also
  \begin{equation*}
  T_n f=1_Y L^n (f1_Y)\ \  \text{ and }\ \  R_n f=1_Y L^n(f1_{\{\phi_\gamma=n\}})
  \end{equation*}
The operator $T_n$ counts all returns from $Y$ to $Y$ at time $n$,
while $R_n$ takes only the first returns at time $n$ into
account. Note that, by definition, $\phi_\gamma=0$ outside of
$Y$, so $R_n$ really counts returns to $Y$. For $z\in \D$, we set
$T(z)=I+\sum T_n z^n$ and $R(z)=\sum R_n z^n$.

\begin{lem}
$T_n$ and $R_n$ are bounded operators on $\boL$, $\norm{T_n}=O(1)$,
$\norm{R_n}=O(m[\phi_\gamma=n])$ and, $\forall z\in \D$, $T(z)=(I-R(z))^{-1}$.
\end{lem}
\begin{proof}
We have $R_n=\sum_{\underline{d}=[d_0]_\gamma,
d_0=[a_0,\ldots,a_{n-1},\gamma]} M_{\underline{d}}$. Thus, Lemma
\ref{lemme_sur_Mdbarre} shows that $\norm{R_n}\leq B(1+\theta)
\sum m[\underline{d}]=(1+\theta) B m[\phi_\gamma=n]$.

In the same way, $T_n=\sum M_{\underline{d}}$ where the sum extends to
all $\underline{d}=[d_0,\ldots,d_{k-1}]_\gamma$ with
$d_i=[\xi_{i0},\ldots,\xi_{in_i},\gamma]$ and $\sum (n_i+1)=n$. Hence,
$\norm{T_n}\leq B(1+\theta)m[Y]$ (the sum is a sum of measures of disjointed
sets included in $Y$, less than $m[Y]$).

Finally, $T_n$ counts all returns to $Y$ while $R_n$ counts only the first
returns. Hence, $T_n=\sum_{i_1+\ldots+i_k=n}R_{i_1}\ldots R_{i_k}$, which gives
the renewal equation.
\end{proof}

\begin{lem}
The operator $R(1):\boL \to \boL$ has a simple isolated
eigenvalue at $1$, the spectral projection being given by $P
f=\frac{1}{m[Y]} \int_Y f \dd m$.
\end{lem}
\begin{proof}
As $R(1)$ counts the first returns to $Y$, it is not hard to check
that $R(1)=L_\gamma$ is the transfer operator associated to $T_\gamma$,
i.e. $R(1)=\sum_{\underline{d}=[d_0]_\gamma} M_{\underline{d}}$. In fact,
$R(1)^n=\sum_{\underline{d}=[d_0,\ldots,d_{n-1}]_\gamma} M_{\underline{d}}$,
  hence Lemma \ref{lemme_sur_Mdbarre} shows that
  \begin{equation}
  \label{lasota_yorke_1}
  \norm{R(1)^n f}_\boL \leq B \theta^n \norm{f}_\boL+B\norm{f}_1
  \end{equation}
The injection $\boL \to L^1(m)$ is compact by the Arzela-Ascoli
theorem. Hence, the Doeblin-Fortet inequality
\eqref{lasota_yorke_1} gives, with the use of Hennion's theorem
(\cite{hennion}), that the essential spectral radius of $R(1)$
acting on $\boL$ is $\leq \theta$. Thus, if $1$ is an eigenvalue
of $R(1)$, it is automatically isolated and of finite
multiplicity.

As $T_\gamma$ preserves the measure $m_\gamma$ (since $T$
preserves $m$), $L_\gamma1=1$ and $PR(1)=R(1)P=P$. By Lemma
\ref{top_mixing_induit_bien} and Proposition
\ref{hypotheses_donnent_exactitude}, $T_\gamma$ is ergodic,
whence there is no other eigenfunction for the eigenvalue $1$.
Finally, there is no nilpotent part for this eigenvalue either,
since $\norm{R(1)^n}$ remains bounded.
\end{proof}

\begin{lem}
$\forall z\in \overline{\D}-\{1\}$, $I-R(z)$ is invertible on $\boL$.
\end{lem}
\begin{proof}
Summing the estimates given by Lemma \ref{lemme_sur_Mdbarre} for
$\underline{d}$ of length $n$ gives that
    \begin{equation}
    \label{df_ineq}
    \norm{R(z)^nf}_{\boL}\leq B|z|^n(\theta^n \norm{f}_\boL+\norm{f}_1)
    \end{equation}
As the injection $\boL \to L^1(m)$ is compact by the Arzela-Ascoli
theorem, the theorem of Hennion (\cite{hennion}) ensures that, $\forall z\in
\overline{\D}$, the
essential spectral radius of $R(z)$ acting on $\boL$ is $\leq
\theta<1$. To obtain the invertibility of $I-R(z)$, it is thus
enough to show that $1$ is not an eigenvalue of $R(z)$. The only
problem is for $|z|=1$ because otherwise, again by
Equation \eqref{df_ineq}, the spectral radius of $R(z)$ is
$\leq |z|<1$ (since $\norm{f}_1\leq \norm{f}_\boL$).
So, let $z=e^{it}$ be fixed, with $0<t<2\pi$.

Suppose that $R(z)f=f$ for some nonzero $f\in \boL$. We will
write, for $u,v\in L^2(m_\gamma)$, $\langle u,v\rangle=\int
\overline{u}v \dd m_\gamma$. Define the operator $W:L^\infty(m_\gamma)\to
L^\infty(m_\gamma)$ by $Wu=e^{-it\phi_\gamma}u\circ T_\gamma$. As
$R(z)v=R(1)(e^{it\phi_\gamma}v)$, this operator $W$ satisfies
    \begin{equation*}
    \langle u,R(z)v\rangle
    =\int \overline{u}\,R(z)v
    =\int \overline{u}\,R(1)(e^{it\phi_\gamma}v)
    =\int \overline{u}\circ T_\gamma\, e^{it\phi_\gamma}v
    =\int \overline{Wu}\cdot v
    =\langle Wu,v\rangle
    \end{equation*}
We show that $f$ is an eigenfunction of $W$ for the eigenvalue $1$:
    \begin{align*}
    \norm{Wf-f}_2^2&
    =\norm{Wf}_2^2-2\Ree\langle Wf,f\rangle+\norm{f}_2^2
    =\norm{Wf}_2^2-2\Ree\langle f,R(z)f\rangle +\norm{f}_2^2
    \\&
    =\norm{Wf}_2^2-2\Ree\langle f,f\rangle+\norm{f}_2^2
    =\norm{Wf}_2^2-\norm{f}_2^2
    \end{align*}
As $T_\gamma$ preserves the measure $m_\gamma$, we have
$\norm{Wf}_2^2=\int |f|^2\circ T_\gamma=\int |f|^2=\norm{f}_2^2$,
which gives $\norm{Wf-f}_2^2=0$. Hence, the function $Wf-f$ is
zero $m_\gamma$-almost everywhere. As $f\in \boL$ and $m_\gamma$
is nonzero on every cylinder, the function $f$ is continuous, thus
$Wf-f=0$ everywhere.

We have a function $f$ such that $e^{-it\phi_\gamma}f\circ
T_\gamma=f$. Taking the modulus, the ergodicity of $T_\gamma$
gives that $|f|$ is constant almost everywhere, hence everywhere
by continuity. As $f\not\equiv 0$, this constant is nonzero, and
we get $e^{-it\phi_\gamma}=f/f\circ T_\gamma$. We can apply
Theorem 3.1. in \cite{aaronson_denker} and obtain that $f$ is
$\delta^*$-measurable, where $\delta^*$ is the smallest partition
such that $\forall d\in \delta$, $T_\gamma d$ is contained in an
atom of $\delta^*$. As every $T_\gamma d$ is a union of sets of
$\gamma$, this implies in particular that $f$ is constant (almost
everywhere, hence everywhere by continuity) on each set of
$\gamma$.

Let $a\in \gamma$. On $[a]$, $f$ is equal to a constant $c$. As
$T$ is topologically mixing, there exists $N$ such that, $\forall
n\geq N$, $[a]\subset T^n[a]$. Let $n\geq N$, and $x\in [a]$ be
such that $T^n x\in [a]$. Let $T^{k_1}x, T^{k_2}x,\ldots,T^{k_p}x$
be the successive returns of $x$ to $Y$, with $k_p=n$. Then
$T^nx=T_\gamma^px$ and $n=\sum_{k=0}^{p-1} \phi_\gamma(T_\gamma^k
x)$. Thus,
  \begin{equation*}
  e^{-itn}=e^{-it \sum_{k=0}^{p-1} \phi_\gamma(T_\gamma^k x)}
  =\frac{f(x)}{f(T_\gamma x)} \frac{f(T_\gamma x)}{f(T_\gamma^2 x)} \cdots
  \frac{f(T_\gamma^{p-1}x)}{f(T_\gamma^px)}
  =\frac{f(x)}{f(T^n x)}=\frac{c}{c}=1
  \end{equation*}
This is true for any $n\geq N$. Taking for example $n=N$ and $N+1$ and
quotienting, we obtain $e^{it}=1$, which is a contradiction.
\end{proof}

\begin{lem}
We have $PR'(1)P=\frac{1}{m[Y]}P$.
\end{lem}
\begin{proof}
Using the
explicit formula for the spectral projection $P$, it is not
difficult to check that $PR_n P=\frac{m[\phi_\gamma=n]}{m[Y]}P$, and
consequently $PR'(1)P=\frac{1}{m[Y]}P$ by the Kac formula
(\cite[Formula 1.5.5]{aaronson:book}). To apply this formula, we have to check
that $T$ is conservative and ergodic, knowing that this is the case for
$T_\gamma$. This can be done for example using \cite[Proposition
1.5.2]{aaronson:book}.
\end{proof}

\begin{proof}[Proof of Theorem \ref{thm_applications_markov_precis}]
The lemmas above show that the hypotheses of Theorem
\ref{theoreme_renouvellement} are satisfied. Consequently, we
get the existence of $E_n \in \Hom(\boL,\boL)$
with $\norm{E_n}=O(F_\beta(n))$ such that $\forall f\in \boL$,
  \begin{equation*}
  1_Y \hat{T}^n f=1_Y\left( \int f \dd m +\sum_{k=n+1}^\infty m[\phi_\gamma >n] \int f
    \dd m +E_n f \right)
  \end{equation*}
Multiplying by an arbitrary $g\in L^\infty(X,\boB,m)$ supported inside $Y$,
we have by the definition of the transfer operator
  \begin{equation*}
  \int f\cdot g\circ T^n \dd m
  =\int f\int g+\sum_{k=n+1}^\infty m[\phi_\gamma>k] \int f\int g+\int g \cdot E_n f \dd m
  \end{equation*}
The absolute value of the last term is bounded by $\norm{g}_\infty
\norm{E_n}_\boL \norm{f}_\boL$, which gives the result.

Finally, if $\int f=0$, we use Theorem \ref{thm_quand_Pf_nul} and
conclude in the same way, the estimates with $F_\beta(n)$ being replaced
by estimates in $O(1/n^\beta)$.
\end{proof}

\subsection{Decay of correlations on the whole space}
Theorem \ref{thm_applications_markov_precis} gives a very sharp
estimate on the decay of correlations when the functions $f$ and
$g$ are supported in $Y$. It is also possible to estimate the
speed of decay for a general $g$, not necessarily supported in
$Y$, although the estimates will be less precise. This kind of
result will be useful in the proof of the Central Limit Theorem.

\begin{thm}
\label{decay_sur_tout_X}
Under the hypotheses of Theorem
\ref{thm_applications_markov_precis}, assume that $f$ is
supported in $Y$ and that $g\in L^\infty(m)$. Then there exists a
constant $C$ (independent of $f$ or $g$) such that
  \begin{equation*}
  \Cor(f,g\circ T^n)\leq \frac{C}{n^{\beta-1}} \norm{f}_\boL
  \norm{g}_\infty
  \end{equation*}
\end{thm}

To obtain this theorem, it is enough to prove that $\|\hat{T}^n f
-\int f\|_1 \leq \frac{C}{n^{\beta-1}}\norm{f}_\boL$.

\begin{lem}
\label{lemme_decay_sur_tout_X}
There exists $C$ such that $\forall
A\in \boB$, $\left|\int_A \hat{T}^n f\dd m -m(A)\int f\right| \leq
\frac{C}{n^{\beta-1}}\norm{f}_\boL $.
\end{lem}
\begin{proof}
In the course of this proof, we shall write $L$ for the transfer
operator acting on functions in $\boL$. Write also $K_0 f=1_Y f$
and, $\forall k\geq 1$, $K_k f=L(1_{\{\phi_\gamma>k\}}f)$: $K_k$
counts the first returns to $Y$ at time $k$, even for points not
starting in $Y$ (contrary to $R_k$). It is then easy to check that
$L^n f=\sum_{k=0}^n K_k T_{n-k} f$ for any $f$ supported in $Y$
(recall that, outside of $Y$, $\phi_\gamma=0$ by definition).

Then, writing $T_n f=\int f +\epsilon_n$ with
$\norm{\epsilon_n}_\boL\leq \frac{C}{n^{\beta-1}}\norm{f}_\boL$,
  \begin{align*}
  \int_A L^n f
  &=\int 1_{Y\cap A} T_n f+\sum_{k=1}^n \int 1_A L^k(1_{\{\phi_\gamma>k\}}T_{n-k}f)
  \\&=\int 1_{Y\cap A} T_n f+\sum_{k=1}^n \int 1_A\circ T^k \cdot 1_{\{\phi_\gamma>k\}}T_{n-k}f
  \\&=\int f \Biggl( \int 1_{Y\cap A} +\sum_{k=1}^n \int 1_A\circ
  T^k\cdot
    1_{\{\phi_\gamma>k\}}\Biggr)
  \\&\hphantom{=}\ + \Biggl(\int 1_{Y\cap A} \epsilon_n +\sum_{k=1}^n \int 1_A\circ T^k
  \cdot 1_{\{\phi_\gamma>k\}}\epsilon_{n-k}\Biggr)
  \\&=I\int f+ II
  \end{align*}
$I$ can be expressed as $m(Y\cap A)+\sum_{k=1}^n m(Y\cap T^{-k}A
-\bigcup_{j=1}^k T^{-j}Y)$. Thus, by Kac's Formula (see
\cite[Lemma 1.5.4]{aaronson:book}), $I=m(A)-\sum_{k=n+1}^\infty
m(Y\cap T^{-k}A-\bigcup_{j=1}^k T^{-j}Y)$. As $m(Y\cap
T^{-k}A-\bigcup_{j=1}^k T^{-j}Y) \leq m[\phi_\gamma>k]\leq
C/k^\beta$ by hypothesis, a summation yields
$I=m(A)+O(1/n^{\beta-1})$.

In the same way, $II\leq \norm{\epsilon_n}_\infty+\sum_{k=1}^n
m[\phi_\gamma>k]\norm{\epsilon_{n-k}}_\infty$: this is a
convolution between sequences respectively in $O(1/n^\beta)$ and
$O(1/n^{\beta-1})$, whence $II=O(1/n^{\beta-1})$ by Lemma
\ref{convole_suites}.
\end{proof}

\begin{proof}[Proof of Theorem \ref{decay_sur_tout_X}]
Lemma \ref{lemme_decay_sur_tout_X} yields that, $\forall A\in
\boB$,
  \begin{equation*}
  \left|\int_A \left( \hat{T}^n f -\int f\right) \dd m \right| \leq
  \frac{C}{n^{\beta-1}}\norm{f}_\boL
  \end{equation*}
Apply this estimate to $A=\{ \hat{T}^n f-\int f \geq 0\}$, then
to $A=\{\hat{T}^n f-\int f <0\}$, and sum to obtain that
  \begin{equation*}
  \int \left|\hat{T}^n f -\int f\right| \dd m \leq
  \frac{2C}{n^{\beta-1}}\norm{f}_\boL
  \end{equation*}
\end{proof}

\subsection{Central Limit Theorem}

\begin{prop}
\label{decay_sur_tout_X_integrale_nulle} Under the hypotheses of
Theorem \ref{thm_applications_markov_precis}, assume that $f$ is
supported in $Y$ and that $\int f=0$. Then there exists a constant
$C$ (independent of $f$) such that $\|\hat{T}^n f\|_1 \leq
\frac{C}{n^\beta}\norm{f}_\boL$.
\end{prop}
\begin{proof}
This is an analogue of Theorem \ref{decay_sur_tout_X} in the case
where $\int f=0$ (which implies that there is a better bound on
$\norm{T_n f}_\boL$, according to Theorem
\ref{thm_applications_markov_precis}). The same proof works
again, and is even easier because the term $I$ in the proof of
Lemma \ref{lemme_decay_sur_tout_X} disappears.
\end{proof}

The following lemma will be useful in the Central Limit Theorem
to precise the regularity of the cocycle in the case of zero
variance.
\begin{lem}
\label{cocycle_regulier} Let $(X,\boB,T,m,\alpha)$ be an
irreducible probability preserving Markov map with the big image
property and for which the distortion $\log g_m$ is Hölderian for
an exponent $\theta<1$. Let $\boL$ denote the space of bounded
functions such that $\sup_{a\in \alpha'} D_a f <+\infty$. If $f\in
\boL$ and $g:X\to \R$ is measurable and satisfies $f\circ
T=g\circ T-g$, then $g\in \boL$.
\end{lem}
\begin{proof}
Denote by $\alpha^n(x)$ the element of the partition
$\bigvee_{i=0}^{n-1}T^{-i}\alpha$ containing $x$. A classical
theorem on continuity points of measurable functions (true on
$[0,1]$, in which $X$ can be canonically imbedded) implies that
  \begin{equation*}
  \text{for almost every }x, \forall \epsilon>0,
  \frac{m\{y\in \alpha^n(x) \tq
  |g(y)-g(x)|>\epsilon\}}{m[\alpha^n(x)]} \to 0 \text{ as }n
  \to \infty
  \end{equation*}
The points that visit infinitely many times every element of the
partition $\alpha$ form also a set of probability $1$. We fix a
point $x_0$ verifying these two properties.

Fix $\epsilon>0$. Let $n_k\to \infty$ be a sequence such that
$T^{n_k}x_0$ visits infinitely often every element of $\alpha$
too, and $\sum \frac{m\{y\in \alpha^{n_k}(x_0) \tq
  |g(y)-g(x_0)|>\epsilon\}}{m[\alpha^{n_k}(x_0)]} <\infty$. For every $k\in \N$, the
control on the distortion implies that
  \begin{multline*}
  \frac{m\{ y\in X \tq \exists y'\in \alpha^{n_k}(x_0), T^{n_k}y'=y,
  |g(y')-g(x_0)|>\epsilon\}}{m[T^{n_k} \alpha^{n_k}(x_0)]}\\
  \asymp
  \frac{m\{y'\in \alpha^{n_k}(x_0) \tq
  |g(y')-g(x_0)|>\epsilon\}}{m[\alpha^{n_k}(x_0)]}
  \end{multline*}
Thus, $\sum_k m\{ y\in X \tq \exists y'\in \alpha^{n_k}(x_0),
T^{n_k}y'=y, |g(y')-g(x_0)|>\epsilon\}<+\infty$. Consequently,
$A_\epsilon:=\{ y\in X\tq \exists K, \forall k\geq K, \text{ if
}y'\in \alpha^{n_k}(x)\text{ is such that }T^{n_k}y'=y,\text{
then }|g(y')-g(x_0)|<\epsilon\}$ is of full measure.

Take $y_1,y_2\in A_\epsilon$ such that $y_1$ and $y_2$ are in the
same element of $\alpha'$. If $d(y_1,y_2)=\theta^n$, it is
possible to write $y_i=[a_0,\ldots,a_{n-1},z_i]$. Take $k$ such
that $T^{n_k}x_0$ is in the same element of $\alpha'$ as $y_1$
and $y_2$. If $k$ is large enough, by definition of $A_\epsilon$,
the preimages $y'_i$ of $y_i$ in $\alpha^{n_k}(x_0)$ satisfy
$|g(y'_i)-g(x_0)|\leq \epsilon$, hence $|g(y'_1)-g(y'_2)|\leq
2\epsilon$. Then
  \begin{align*}
  |g(y_1)-g(y_2)|&=
  |g\circ T^{n_k}(y'_1)-g\circ T^{n_k}(y'_2)|
  \\&
  \leq \sum_{i=1}^{n_k} |f\circ T^i(y'_1)-f\circ T^i(y'_2)|
  +|g(y'_1)-g(y'_2)|
  \\
  &\leq \sum_{i=1}^{n_k} \norm{f}_\boL \theta^{n_k+n-i} +2\epsilon
  \leq \frac{\norm{f}_\boL}{1-\theta}d_\theta(y_1,y_2)+2\epsilon
  \end{align*}

Finally, for $y_1,y_2\in A=\bigcap A_\epsilon$ of full measure,
$|g(y_1)-g(y_2)|\leq \frac{\norm{f}_\boL}{1-\theta}d(y_1,y_2)$.
Hence, there exists a unique version of the function $g$ which is
Lipschitzian on every set of $\alpha'$, which we will still
denote by $g$.

To see that $g\in \boL$, there remains to prove that $g$ is
bounded. Let $\eta>0$ be such that $\forall a\in \alpha,
m[Ta]>\eta$. There exists $a_1,\ldots,a_N \in \alpha$ a finite
number of partition sets such that $\sum_{i=1}^N m[a_i]>1-\eta$.
Thus, $\forall a\in \alpha$, $Ta$ contains one of the sets $a_i$.
On each of these sets, $g$ is Lipschitzian, hence bounded by a
constant $C_i$. If $x\in [a]$ has its image in $a_i$, then
$|g(x)|=|g\circ T(x)-f\circ T(x)| \leq C_i+\norm{f}_\infty\leq
\max_i C_i+\norm{f}_\infty=:C$. Finally, for $y\in [a]$,
$|g(y)|\leq |g(x)-g(y)|+|g(y)|\leq
\frac{\norm{f}_\boL}{1-\theta}+C$.
\end{proof}

\begin{thm}
\label{markov_CLT}
Under the hypotheses of Theorem
\ref{thm_applications_markov_precis}, if $f\in \boL$ is supported
in $Y$ and $\int f=0$, then the sequence $\frac{1}{\sqrt{n}}
\sum_{k=0}^{n-1} f\circ T^k$ converges in law to a Gaussian
random variable of zero mean and finite variance $\sigma$, with
  \begin{equation*}
  \sigma^2=-\int f^2\dd m+2\sum_{n=0}^\infty \int f\cdot f\circ
  T^n \dd m
  \end{equation*}
Moreover, $\sigma=0$ if and only if there exists a measurable
function $g$ such that $f\circ T=g\circ T-g$. Such a function $g$
automatically satisfies $g_{|Y}\in \boL$ and $\forall x\in Y,
\forall n<\phi_\gamma(x), g(T^nx)=g(x)$.
\end{thm}

We will use an abstract result due to Liverani \cite[Theorem
1.1]{liverani:CLT} inspired by Kipnis-Varadhan to obtain this
Central Limit Theorem. We recall for the convenience of the
reader the version of this theorem that will be useful in our
setting.
\begin{thm}
\label{abstract_CLT}
Let $(X,\boB,T,m)$ be a non-singular
probability preserving dynamical system. Let also $f\in
L^\infty(X)$, $\int f=0$ be such that
\begin{enumerate}
\item $\sum_{n=0}^\infty \left|\int f\cdot f\circ
T^n\right|<\infty$.
\item The series $\sum_{n=0}^\infty \hat{T}^n f$ converges
absolutely in $L^1$.
\end{enumerate}
Then the sequence $\frac{1}{\sqrt{n}} \sum_{k=0}^{n-1} f\circ
T^k$ converges in law to a Gaussian random variable of zero mean
and finite variance $\sigma$, with
  \begin{equation*}
  \sigma^2=-\int f^2\dd m+2\sum_{n=0}^\infty \int f\cdot f\circ
  T^n \dd m
  \end{equation*}
Moreover, $\sigma=0$ if and only if there exists a measurable
function $g$ such that $f\circ T=g\circ T-g$.
\end{thm}

\begin{proof}[Proof of Theorem \ref{markov_CLT}]
It is enough to show that the hypotheses of Theorem
\ref{abstract_CLT} are verified. As we have formulated this
theorem, the first hypothesis is in fact a consequence of the
second one, since
  \begin{equation*}
  \left|\int f \cdot f\circ T^n \dd m\right|
  =\left|\int \hat{T}^n f \cdot f \dd m \right|
  \leq \|\hat{T}^n f\|_1 \norm{f}_\infty
  \end{equation*}
Consequently, there remains only to check that $\sum \|\hat{T}^n
f\|_1 <+\infty$. By Proposition
\ref{decay_sur_tout_X_integrale_nulle}, $\|\hat{T}^n
f\|_1=O(1/n^\beta)$ with $\beta>1$, thus the series is summable.

To obtain the regularity results on $g$ when $\sigma=0$, we use
the fact that $f=0$ outside of $Y$. As $f\circ T=g\circ T-g$,
this implies that $g(x)=g\circ T(x)$ when $T(x)\not\in Y$. In
particular, $\forall x\in Y, \forall n<\phi_\gamma(x), g(x)=g(T^n
x)$. Using once more the cocycle relation gives that $f\circ
T_\gamma(x)=g\circ T_\gamma(x)-g(x)$. Thus, Lemma
\ref{cocycle_regulier} applied to $(Y,T_\gamma)$ shows that
$g_{|Y}\in \boL$.
\end{proof}

\section{Applications to specific maps}
\label{section_applications}
\subsection{The Liverani-Saussol-Vaienti map}
The Liverani-Saussol-Vaienti map is the map $T:[0,1]\to [0,1]$ defined by
   \begin{equation*}
   T(x)=\left\{ \begin{array}{cl}
   x(1+2^\alpha x^\alpha) &\text{if }0\leq x\leq 1/2
   \\
   2x-1 &\text{if }1/2<x\leq 1
   \end{array}\right.
   \end{equation*}

It is shown in \cite{liverani_saussol_vaienti} that, when
$0<\alpha<1$, $T$ admits an integrable invariant density $h$ which
is Lipschitz outside of any neighborhood of $0$.

\begin{cor}
\label{cor_LSV}
If $\alpha\in (0,1)$, $f$ is Lipschitz, $g$ is
bounded measurable, $\int f,\int g\not=0$ and $f,g=0$ in a
neighborhood of $0$, then $\Cor(f,g\circ T^n)\sim
\frac{1}{4}h(\frac{1}{2}) \alpha^{-1/\alpha}
(\frac{1}{\alpha}-1)^{-1}n^{1-1/\alpha}\int f\int g$ with respect
to the invariant probability measure.

Moreover, if $\int f=0$ (and $f,g$ are still zero in a
neighborhood of $0$, $f$ Lipschitzian), then $\Cor(f,g\circ
T^n)=O(1/n^{1/\alpha})$. Consequently, $f$ satisfies a Central
Limit Theorem.
\end{cor}
\begin{proof}
If $x_0=1/2$ and $x_{i+1}=T^{-1}(x_i)\cap [0,1/2]$, the partition
$\alpha=\{(x_{i+1},x_i]\}\cup (1/2,1]$ is a Markov partition for
$T$, which makes it possible to apply the results of the previous
section to this map. The distortion of the induced map on
$(1/2,1]$ is locally Hölder continuous for the density $h$,
whence Theorem \ref{thm_applications_markov_precis} applies and
gives a precise asymptotic on the speed of decay of correlations
for functions supported in $(1/2,1]$, which can be calculated
precisely (see \cite{sarig:decay}).

As the distortion from $(x_{i+1},x_i]$ to $(x_i,x_{i-1}]$ is bounded,
it is not hard to check that the induced map on $\gamma=\{(x_{i+1},x_i]
\tq i<N\} \cup \{(1/2,1]\}$ has still a Hölder continuous
distortion for any $N$.
Thus, Theorem \ref{thm_applications_markov_precis} gives
also estimates on the decay of correlations of functions supported in
$(x_N,1]$. More precisely, for functions $f\in \boL$ and $g\in L^\infty$
supported in $(x_N,1]$,
  \begin{equation}
  \label{liv_1}
  \Cor(f,g\circ T^n)\sim \left(
  \sum_{k=n+1}^{\infty}m[\phi_\gamma>k]\right)\int f \int
  g
  \end{equation}
For functions supported in $(1/2,1]$, Sarig has shown, estimating
$m[\phi_{(1/2,1]}]>n$, that
  \begin{equation}
  \label{liv_2}
  \Cor(f,g\circ T^n)\sim \frac{1}{4}h\left(\frac{1}{2}\right) \alpha^{-1/\alpha}
  \left(\frac{1}{\alpha}-1\right)^{-1}n^{1-1/\alpha}\int f\int g
  \end{equation}
The estimate \eqref{liv_1} can be applied in particular to functions
supported in $(1/2,1]$, which gives, after comparing with
\eqref{liv_2}, that $\sum_{k=n+1}^{\infty}m[\phi_\gamma>k] \sim \frac{1}{4}h(\frac{1}{2}) \alpha^{-1/\alpha}
  (\frac{1}{\alpha}-1)^{-1}n^{1-1/\alpha}$. This proves the corollary.
\end{proof}

\subsection{LS Young towers}
\label{section_LSY}
A \emph{LS Young tower} is a non-singular
conservative transformation $(\Delta,\boB,m,F)$ with a generating
partition $\{\Delta_{l,i} \tq i\in \N, l=0,\ldots,R_i-1\}$ with
the following properties:
\begin{enumerate}
  \item $\forall l,i$ the measure of $\Delta_{l,i}$ is positive and finite.
    Moreover, if $\Delta_l=\bigcup \Delta_{l,i}$,
    $m(\Delta_0)<\infty$.
  \item If $l+1<R_i$, $F:\Delta_{l,i}\to \Delta_{l+1,i}$ is a measurable
    bijection and $F_* m_{|\Delta_{l,i}}=m_{|\Delta_{l+1,i}}$.
  \item If $l+1=R_i$, $F:\Delta_{l,i}\to \Delta_0$ is a measurable bijection.
  \item Let $R:\Delta_0\to \N$ be the function $R_{|\Delta_{0,i}}=R_i$, and set
    $g=\frac{\de m_{|\Delta_0}}{\de m_{|\Delta_0}\circ F^R}$. $g$ has a version
    for which $\exists C>0, \theta \in (0,1)$ such that $\forall i$ and $\forall x,y\in
    \Delta_{0,i}$,
    \begin{equation*}
    \left|\frac{g(x)}{g(y)}-1\right|\leq C \theta^{s(F^Rx,F^Ry)}
    \end{equation*}
    where $s(x,y)=\min\{n \tq (F^R)^nx, (F^R)^ny \text{ lie in different
    }\Delta_{0,j}\}$.
\end{enumerate}
The fourth condition corresponds exactly to saying that the induced
map on the base $\Delta_0$ of the tower has a distortion which is locally
Hölder continuous.

Henceforth, we assume for simplicity that $\int R \dd m<+\infty$ and that
$m$ is an $F$-invariant probability, which is possible because $m$ has
an integrable invariant density $h$ such that $c_0^{-1}\leq h \leq c_0$
(see \cite[Theorem 1]{lsyoung:recurrence}).

Set $C_\theta(\Delta)=\{f:\Delta \to \C \tq \exists C \forall x,y\in \Delta, |f(x)-f(y)|\leq
C\theta^{s(x,y)}\}$: this is the space of locally Hölder continuous
functions ($s$ has been extended to all pairs $x,y\in \Delta$ by setting
$s(x,y)=0$ if $x,y$ are not in the same $\Delta_{l,i}$ and, for $x,y\in
\Delta_{l,i}$, $s(x,y)=s(x',y')$ where $x',y'$ are the corresponding
points in $\Delta_{0,i}$).

\begin{cor}
\label{cor_LSY}
Let $(\Delta,\boB,m,F)$ be a probability preserving LS Young tower with
$\gcd\{R_i\}=1$ and $m[R>n]=O(1/n^\beta)$ where $\beta>1$. If
$f\in C_\theta(\Delta)$, $g\in L^{\infty}$ are supported inside
$\bigcup_0^{N-1}\Delta_l$ for some $N$, then $\Cor(f,g\circ
F^n)=\sum_{k>n} m[R>k]\int f\int g+O(F_\beta(n))$.

Moreover, if $\int f=0$, then $\Cor(f,g\circ F^n)=O(1/n^\beta)$. Thus,
$f$ satisfies a Central Limit Theorem.
\end{cor}
\begin{proof}
For the partition $\{\Delta_{l,i}\}$, $F$ does not have the big image
  property. However, it is still a Markov map for the partition
  $\{\Delta_l\}$ composed of the points at different heights. If
  $\gamma=\{\Delta_l \tq l<N\}$ for some $N$, then $\gamma$ is finite, whence the
  induced map $T_\gamma$ has the big image property.

For the induced map, the partition $\delta$ is constructed as
follows: at each height $0<l<N-1$, cut $\Delta_l$ in two pieces
$\Delta_l\cap F^{-1}\Delta_0$ and $\Delta_l-F^{-1}\Delta_0$.
$\Delta_0$ remains intact, and $\Delta_{N-1}$ is cut into all the
small pieces $\Delta_{N-1,i}$. With this explicit partition, it
is not hard to check that the induced map has
$\theta^{1/N}$-locally Hölder continuous distortion.

Thus, Theorem \ref{thm_applications_markov_precis} applies and gives
an estimate
  \begin{equation}
  \label{eqn:equiv_dans_tours}
  \Cor(f,g\circ F^n)=\sum_{k>n} m[\phi_\gamma>k]\int f\int g+O(F_\beta(n))
  \end{equation}
To finish the proof of the theorem, we have to show that
$\sum_{k>n}m[\phi_\gamma>k]=\sum_{k>n}m[R>k]+ O(F_\beta(n))$. If
$f$ and $g$ are supported in $\Delta_0$ and of nonzero integral,
Estimate \eqref{eqn:equiv_dans_tours} applies. Moreover, the
estimate for $N=1$ applies also. Equating these two estimates of
$\Cor(f,g\circ F^n)$, we get the result.
\end{proof}

\bibliography{biblio}
\bibliographystyle{alpha}
\end{document}